\documentclass{amsart}
\usepackage{latexsym}
\usepackage{amssymb}
\usepackage{amsmath}

\begin{document}


\newtheorem{theorem}{Theorem}
\newtheorem{problem}{Problem}
\newtheorem{definition}{Definition}
\newtheorem{lemma}{Lemma}
\newtheorem{proposition}{Proposition}
\newtheorem{corollary}{Corollary}
\newtheorem{example}{Example}
\newtheorem{conjecture}{Conjecture}
\newtheorem{algorithm}{Algorithm}
\newtheorem{exercise}{Exercise}
\newtheorem{remarkk}{Remark}

\newcommand{\be}{\begin{equation}}
\newcommand{\ee}{\end{equation}}
\newcommand{\bea}{\begin{eqnarray}}
\newcommand{\eea}{\end{eqnarray}}
\newcommand{\beq}[1]{\begin{equation}\label{#1}}
\newcommand{\eeq}{\end{equation}}
\newcommand{\beqn}[1]{\begin{eqnarray}\label{#1}}
\newcommand{\eeqn}{\end{eqnarray}}
\newcommand{\beaa}{\begin{eqnarray*}}
\newcommand{\eeaa}{\end{eqnarray*}}
\newcommand{\req}[1]{(\ref{#1})}

\newcommand{\lip}{\langle}
\newcommand{\rip}{\rangle}

\newcommand{\uu}{\underline}
\newcommand{\oo}{\overline}
\newcommand{\La}{\Lambda}
\newcommand{\la}{\lambda}
\newcommand{\eps}{\varepsilon}
\newcommand{\om}{\omega}
\newcommand{\Om}{\Omega}
\newcommand{\ga}{\gamma}
\newcommand{\rrr}{{\Bigr)}}
\newcommand{\qqq}{{\Bigl\|}}

\newcommand{\dint}{\displaystyle\int}
\newcommand{\dsum}{\displaystyle\sum}
\newcommand{\dfr}{\displaystyle\frac}
\newcommand{\bige}{\mbox{\Large\it e}}
\newcommand{\integers}{{\Bbb Z}}
\newcommand{\rationals}{{\Bbb Q}}
\newcommand{\reals}{{\rm I\!R}}
\newcommand{\realsd}{\reals^d}
\newcommand{\realsn}{\reals^n}
\newcommand{\NN}{{\rm I\!N}}

\newcommand{\degree}{{\scriptscriptstyle \circ }}
\newcommand{\dfn}{\stackrel{\triangle}{=}}
\def\complex{\mathop{\raise .45ex\hbox{${\bf\scriptstyle{|}}$}
     \kern -0.40em {\rm \textstyle{C}}}\nolimits}
\def\hilbert{\mathop{\raise .21ex\hbox{$\bigcirc$}}\kern -1.005em {\rm\textstyle{H}}} 
\newcommand{\RAISE}{{\:\raisebox{.6ex}{$\scriptstyle{>}$}\raisebox{-.3ex}
           {$\scriptstyle{\!\!\!\!\!<}\:$}}} 

\newcommand{\hh}{{\:\raisebox{1.8ex}{$\scriptstyle{\degree}$}\raisebox{.0ex}
           {$\textstyle{\!\!\!\! H}$}}}

\newcommand{\OO}{\won}
\newcommand{\calA}{{\mathcal A}}
\newcommand{\BB}{{\mathcal B}}
\newcommand{\calC}{{\cal C}}
\newcommand{\calD}{{\cal D}}
\newcommand{\calE}{{\cal E}}
\newcommand{\calF}{{\mathcal F}}
\newcommand{\calG}{{\cal G}}
\newcommand{\calH}{{\cal H}}
\newcommand{\calK}{{\cal K}}
\newcommand{\calL}{{\mathcal L}}
\newcommand{\calM}{{\cal M}}
\newcommand{\calO}{{\cal O}}
\newcommand{\calP}{{\cal P}}
\newcommand{\calU}{{\mathcal U}}
\newcommand{\calX}{{\cal X}}
\newcommand{\calXX}{{\cal X\mbox{\raisebox{.3ex}{$\!\!\!\!\!-$}}}}
\newcommand{\calXXX}{{\cal X\!\!\!\!\!-}}
\newcommand{\gi}{{\raisebox{.0ex}{$\scriptscriptstyle{\cal X}$}
\raisebox{.1ex} {$\scriptstyle{\!\!\!\!-}\:$}}}
\newcommand{\intsim}{\int_0^1\!\!\!\!\!\!\!\!\!\sim}
\newcommand{\intsimt}{\int_0^t\!\!\!\!\!\!\!\!\!\sim}
\newcommand{\pp}{{\partial}}
\newcommand{\al}{{\alpha}}
\newcommand{\sB}{{\cal B}}
\newcommand{\sL}{{\cal L}}
\newcommand{\sF}{{\cal F}}
\newcommand{\sE}{{\cal E}}
\newcommand{\sX}{{\cal X}}
\newcommand{\R}{{\rm I\!R}}
\renewcommand{\L}{{\rm I\!L}}
\newcommand{\vp}{\varphi}
\newcommand{\N}{{\rm I\!N}}
\def\ooo{\lip}
\def\ccc{\rip}
\newcommand{\ot}{\hat\otimes}
\newcommand{\rP}{{\Bbb P}}
\newcommand{\bfcdot}{{\mbox{\boldmath$\cdot$}}}

\renewcommand{\varrho}{{\ell}}
\newcommand{\dett}{{\textstyle{\det_2}}}
\newcommand{\sign}{{\mbox{\rm sign}}}
\newcommand{\TE}{{\rm TE}}
\newcommand{\TA}{{\rm TA}}
\newcommand{\won}{{\mbox{\bf 1}}}
\newcommand{\Lebn}{{\rm Leb}_n}
\newcommand{\Prob}{{\rm Prob\,}}
\newcommand{\sinc}{{\rm sinc\,}}
\newcommand{\ctg}{{\rm ctg\,}}
\newcommand{\loc}{{\rm loc}}
\newcommand{\trace}{{\,\,\rm trace\,\,}}
\newcommand{\Dom}{{\rm Dom}}
\newcommand{\ifff}{\mbox{\ if and only if\ }}
\newcommand{\nproof}{\noindent {\bf Proof:\ }}
\newcommand{\remark}{\noindent {\bf Remark:\ }}
\newcommand{\remarks}{\noindent {\bf Remarks:\ }}
\newcommand{\note}{\noindent {\bf Note:\ }}

\newcommand{\boldx}{{\bf x}}
\newcommand{\boldX}{{\bf X}}
\newcommand{\boldy}{{\bf y}}
\newcommand{\boldR}{{\bf R}}
\newcommand{\uux}{\uu{x}}
\newcommand{\uuY}{\uu{Y}}

\newcommand{\limn}{\lim_{n \rightarrow \infty}}
\newcommand{\limN}{\lim_{N \rightarrow \infty}}
\newcommand{\limr}{\lim_{r \rightarrow \infty}}
\newcommand{\limd}{\lim_{\delta \rightarrow \infty}}
\newcommand{\limM}{\lim_{M \rightarrow \infty}}
\newcommand{\limsupn}{\limsup_{n \rightarrow \infty}}

\newcommand{\ra}{ \rightarrow }

\newcommand{\ARROW}[1]
  {\begin{array}[t]{c}  \longrightarrow \\[-0.2cm] \textstyle{#1} \end{array} }

\newcommand{\AR}
 {\begin{array}[t]{c}
  \longrightarrow \\[-0.3cm]
  \scriptstyle {n\rightarrow \infty}
  \end{array}}

\newcommand{\pile}[2]
  {\left( \begin{array}{c}  {#1}\\[-0.2cm] {#2} \end{array} \right) }

\newcommand{\floor}[1]{\left\lfloor #1 \right\rfloor}

\newcommand{\mmbox}[1]{\mbox{\scriptsize{#1}}}

\newcommand{\ffrac}[2]
  {\left( \frac{#1}{#2} \right)}

\newcommand{\one}{\frac{1}{n}\:}
\newcommand{\half}{\frac{1}{2}\:}

\def\le{\leq}
\def\ge{\geq}
\def\lt{<}
\def\gt{>}

\def\squarebox#1{\hbox to #1{\hfill\vbox to #1{\vfill}}}
\newcommand{\nqed}{\hspace*{\fill}
           \vbox{\hrule\hbox{\vrule\squarebox{.667em}\vrule}\hrule}\bigskip}

\newcommand{\no}{\noindent}
\newcommand{\E}{\mathcal{E}}
\newcommand{\EE}{\mathbb{E}}
\newcommand{\RR}{\mathbb{R}}
\newcommand{\D}{\mathcal{D}}
\newcommand{\DD}{\widetilde{\mathcal{D}}}
\newcommand{\DDD}{\widetilde{\widetilde{\mathcal{D}}}}
\newcommand{\LL}{\mathcal{L}}
\newcommand{\F}{\mathcal{F}}
\newcommand{\B}{\mathcal{B}}
\newcommand{\M}{\mathcal{M}}
\newcommand{\W}{\mathbb{W}}
\newcommand{\lo}{lp}
\newcommand{\grandX}{\mathbb{X}}
\newcommand{\loigrandX}{\mu^{\mathbb{X}}}
\newcommand{\mup}{\mu_*}
\newcommand{\betap}{\beta_*}
\newcommand{\pa}{pa}
\newcommand{\QQ}{\mathbb{Q}}


\title{A GENERAL FRAMEWORK FOR VARIATIONAL CALCULUS ON WIENER SPACE}

\author{ K\'evin HARTMANN}
\maketitle
\noindent
{\bf Abstract:}{\small{ We provide a framework to derive a variational formulation for $-\log\EE_\nu\left[e^{-f}\right]$ for a large class of measures $\nu$. We use a family of perturbations of the identity $(W^u)$ whose invertibility we characterize thanks to entropy. This yields results of strong existence for various stochastic differential equations. We also discuss the attainability of the infimum in the variational formulation and we derive a Pr\'ekopa-Leindler theorem for the measure $\nu$.
}}\\

\vspace{0.5cm}

\no Keywords: Wiener space, variational formulation, entropy, invertibility, Brownian bridge, loop measure, diffusing particles, stochastic differential equations, Pr\'ekopa-Leindler theorem
\tableofcontents

\section{\bf{Introduction}}

Denote $\W$ the space of continuous functions from $[0,1]$ to $\RR^n$ and H the associated canonical Cameron-Martin space of elements of $\W$ which admit a density in $L^2$. Also denote $\mu$ the Wiener measure, W the coordinate process, and $(\F_t)$ the canonical filtration of W completed with respect to $\mu$. W is a Brownian motion under $\mu$. Set f a bounded from above measurable function from $\W$ to $\RR$. In \cite{du}, Dupuis and Ellis prove that
\bea -\log\EE_\mu\left[e^{-f}\right]=\inf_\theta\left(\EE_\theta\left[f\right]+H(\theta|\mu)\right)\label{2rv0}\eea
\no where the infimum is taken over the probability measures $\theta$ on $\W$ which are absolutely continuous with respect to $\mu$ and the relative entropy $H(\theta|\mu)$ is equal to $\EE_\mu\left[\frac{d\theta}{d\mu}\log\frac{d\theta}{d\mu}\right]$. In \cite{bd}, Bou\'e and Dupuis use it to derive the variational formulation
\bea -\log\EE_\mu\left[e^{-f}\right]=\inf_u\EE_\mu\left[f\circ (W+u)+\frac{1}{2}\int_0^1|\dot{u}(s)|^2ds\right]\label{2rv}\eea
\no where the infimum is taken over $L^2$ functions from $\W$ to H whose density is adapted to $(\F_t)$. This variational formulation is useful to derive large deviation asymptotics as Laplace principles for small noise diffusions for instance. This result was later extended by Budhiraja and Dupuis to Hilbert-space-valued Brownian motions in \cite{bud}, and then by Zhang to abstract Wiener spaces in \cite{zh}, using the framework developed by \"Ust\"unel and Zakai in \cite{ust1}.\newline
\no The Pr\'ekopa-Leindler theorem first formulation was given by Pr\'ekopa in \cite{pre1} and arose in stochastic programming where a lot of non-linear optimization problems require concavity. In \cite{ha}, Huu Hariya uses the variational formulation to retrieve a Pr\'ekopa-Leindler theorem for Wiener space, similar to the formulation of \"Ust\"unel and Feyel in \cite{ust3} with log-concave measures. Other functional inequalities can be derived from \ref{2rv}, see for instance Lehec in \cite{le1}.\newline
\no The bounded from above hypothesis in \ref{2rv} was weakened significantly by \"Ust\"unel in \cite{art}, it was replaced with the condition
$$\EE_\mu\left[fe^{-f}\right]<\infty$$
\no and the existence of conjugate integers p and q such that
$$f\in L^p(\mu),e^{-f}\in L^q(\mu)$$
\no These relaxed hypothesis pave the way to new applications. The possibility of using unbounded functions is primordial in Dabrowski's application of \ref{2rv} to free entropy in \cite{da}.\newline
\no \"Ust\"unel's approach is routed in the study of the perturbations of the identity of $\W$, which is the coordinate process, and their invertibility. The question of the invertibility of an adapted perturbation of the identity is linked to the representability of measures and was put to light by the celebrated example of Tsirelson \cite{tsi}. \"Ust\"unel proved that if $u\in L^2(\mu,H)$ has an adapted density, $I_\W+u$ is $\mu$-a.s. invertible if and only if
$$H((I_\W+u)\mu|\mu)=\frac{1}{2}\EE_\mu\left[|u|_H^2\right]$$
\no If u satisfies the hypothesis of Girsanov theorem, this gives a criteria of existence of strong solutions to some stochastic differential equations. Indeed, \"Ust\"unel proves in \cite{ust2} that if such a $I_\W+u$ in invertible, its inverse V is a strong solution to the stochastic differential equation
$$dV(t)=-\dot{u}(t)\circ V dt+dW(t)$$
\no To prove \ref{2rv} with the integrability conditions specified above, \"Ust\"unel uses the fact that H-$C^1$ shifts, meaning shifts that are a.s. Fr\'echet-differentiable on H with a a.s. continuous on H Fr\'echet derivative, are a.s. invertible, and that shifts can be approached with H-$C^1$ shifts using the Ornstein-Uhlenbeck semigroup.\newline
\no In \cite{a1} we give a variational formulation similar as \ref{2rv} for diffusions which are solutions of stochastic differential equations, while lowering the integrability hypothesis on f. This paper also focus on the invertibility of certain perturbations of the identity that are not affine shifts since the measure considered is not the same. However, contrary to what \"Ust\"unel does in \cite{art}, we do not approach f to derive invertible shift from those approached functions. We write $\frac{e^{-f}}{\EE\left[e^{-f}\right]}$ as the Wick exponential of some v, and then approach v to obtain shifts that generate invertible perturbations of the identity. Our method relies on the fact that we have a Girsanov-like change of variable formula with the perturbations of the identity, with relatiowith respect to a particular Brownian motion. It does not use any tool that is specific to Gaussian measures.\newline
\no Two questions arise from this: can this method be applied to other measures, and can invertibility results be linked to the existence of strong solutions for some stochastic differential equations?\newline
\no This paper presents a general framework to be able to similarly derive a variational formulation for $-\log\EE_\nu[e^{-f}]$ for a large class of measures $\nu$. We give a set of conditions so that a set of processes $(W^u)$ can act as perturbations of W and allow a Girsanov-like change of variable with respect to a Brownian motion $\beta$. At first we want to have a minimal setting to be able to compute the variational formula, and we just consider the u whose density is a.s. bounded and we prove that
$$-\log\EE_\nu\left[e^{-f}\right]=\inf_u\EE_\nu\left[f\circ W^u+\frac{1}{2}|u|_H^2\right]$$
\no where the infimum is just taken over the u with a.s. bounded density, thus providing a clearer description of the infimum. The integrability hypothesis over f remain the same as in the case of a diffusion. In a second time, we study the possibilities to expand the domain over which the infimum is taken, for both variational calculus results, mainly concerning the attainability of the infimum, and invertibility results. Indeed, we prove that once again, $W^u$ is invertible if and only if
$$H(W^u\nu|\nu)=\frac{1}{2}\EE_\nu\left[|u|_H^2\right]$$
\no and in case of invertibility its inverse is of the form $W^v$. Furthermore, the invertibility of $W^u$ can be related to the existence of strong solutions of stochastic differential equations in certain cases. If $W^u$ can be written $W+w^u$, with $w^u\in L^0(\nu,H)$ having an adapted density, and is invertible, its inverse $W^v$ verify
$$dW^v(t)=-\dot{\overbrace{w^u}}(t)\circ W^vdt+W(t)$$
\no We also prove a Pr\'ekopa-Leindler theorem on $\W$ for the measure $\nu$, however the convexity hypothesis seem very restrictive.\newline
\no We apply this framework to various examples. First we retrieve the case of the image measure of a diffusion of \cite{a1}, and then we study the case of the image measure of a Brownian bridge, a loop measure, and finally the image measure of a set of diffusing particles. The behaving of diffusing particles satisfying a differential stochastic system was studied by C\'epa and Lepingle in \cite{cl}, and Rogers and Shi in \cite{shi}. We focus on the case where the stochastic differential system the particles $(Z_1,...,Z_n)$ verify is of the form
\beaa Z_1(t)&=&z_1(0)+\sigma B_1(t)+b\int_0^t Z_1(s)ds+ct+\gamma\sum_{j\in\{1,...,n\}\backslash\{1\}}\int_0^t\frac{ds}{Z_1(s)-Z_j(s)}\\
\vdots&&\\
Z_n(t)&=&z_n(0)+\sigma B_n(t)+b \int_0^tZ_n(s)ds+ct+\gamma\sum_{j\in\{1,...,n\}\backslash\{n\}}\int_0^t\frac{ds}{Z_n(s)-Z_j(s)}\eeaa
\no where $(B_1,...,B_n)$ is a $\RR^n$-valued Brownian motion and $\sigma^2\leq 2\gamma$, which guarantee there is no collision.

\section{\bf{Framework}\label{fr}}

\no Set $n\in\NN^*$, we denote $\W=C([0,1],\RR^{n})$ the canonical
Wiener space, $H=\left\{\int_0^.\dot{h}(s)ds, \dot{h}\in
  L^2([0,1])\right\}$ the associated Cameron-Martin space and $(W(t))$
is the coordinate process.\newline
\no We assume that $\W$ is equipped with a probability measure $\nu$. The filtration of a process m will be denoted
$\left(\F^m_t\right)$, the filtration of W will be simply denoted
$\left(\F_t\right)$. Except if stated otherwise, every filtration
considered is completed with respect to $\nu$. We denote, for $p\in\RR_+$,
\beaa L^p_a(\nu,H)=\left\{u\in L^p(\nu,H),\dot{u}\;is\;(\F_t)-adapted\right\}\eeaa
\no and
$$\D=\left\{u\in L^0_a(\nu,H), \dot{u}\;is\;d\nu\times
  dt-a.s.\;bounded\right\}$$
\no
\no If m is a martingale and v has a density whose stochastic
integral with respect to m is well defined we will denote
$$\delta_m v =\int_0^1\dot{v}(s)dm(s)$$
\no We also denote the Wick exponential as follow
$$\rho(\delta_m
v)=\exp\left(\int_0^1\dot{v}(s)dm(s)-\frac{1}{2}\int_0^1\left|\dot{v}(s)\right|^2d\langle
  m\rangle(s)\right)$$
\no and for $p\geq 0$ we denote
$$G_p(\nu,m)=\left\{u\in L_a^p(\nu,H), \EE_\nu\left[\rho(-\delta_m u)\right]=1\right\}$$
\no We assume there exists a family of processes
$\left(W^u\right)_{u\in\D}$ and a $\nu$-Brownian motion $\beta$ which verify the following conditions:\newline
\no (i) $\beta$ is a $\nu$-Brownian motion whose canonical filtration
is identical to the canonical filtration of $(W(t))$\newline
\no (ii) $W^0=W$\newline
\no (iii) For every $u\in\D$, the law of $W^u$
under $\tilde{\nu}^u$ is the same as the law of W under $\nu$, where
$\tilde{\nu}^u$ is defined by
$\frac{d\tilde{\nu}^u}{d\nu}=\rho(-\delta_\beta u)$\newline
\no (iv) For every $u\in\D$,
$$\beta\circ W^u=\beta+u$$\newline
\no (v) For every $u,v\in\D$,
$$W^u\circ W^v=W^{v+u\circ W^v}\;\nu-a.s.$$

\begin{remark}
\no Clearly $\D\subset L^\infty_a(\nu,H)$, so if $u\in\D$,
$\EE_\nu\left[\rho(-\delta_\beta u)\right]=1$ and $\tilde{\nu}^u$ which
was defined in condition (iii) is indeed a probability measure.
\end{remark}

\no Condition (iii) can be written as follow:

\begin{proposition}
\label{gir}
\no Set $u\in \D$, for every bounded measurable function f, we have:
$$\EE_\nu\left[f\right]=\EE_\nu\left[f\circ W^u\rho(-\delta_\beta
  u)\right]$$
\end{proposition}

\no Next proposition ensures that the compositions written in (iv) and
(v) are well
defined.

\begin{proposition}
Set $u\in \D$, we have
$$W^u\nu\sim\nu$$
\end{proposition}

\nproof Set $f\in C_b(\W)$ bounded and measurable, we have, using
proposition \ref{gir}
\beaa
\EE_{W^u\tilde{\nu}^u}\left[f\right]&=&\EE_{\tilde{\nu}^u}\left[f\circ W^u\right]\\
&=&\EE_{\nu}\left[f\circ
  W^u\rho\left(-\delta_\beta u\right)\right]\\
&=&\EE_{\nu}\left[f\right]\eeaa
\no so $W^u\tilde{\nu}^u=\nu$.\newline
\no Since $\tilde{\nu}^u\sim\nu$, we have $W^u\tilde{\nu}^u\sim W^u\nu$
which conclude the proof.\nqed

\section{\bf{Invertibility results}}
\subsection{\bf{First results}}
\begin{definition}
A measurable map $U:\W\rightarrow \W$ is said to be
$\nu$-a.s. left-invertible if and only if $U\nu\ll\nu$ and there exists a measurable map
$V:\W\rightarrow \W$ such that $V\circ U=I_\W$ $\nu$-a.s.\newline
\no A measurable map $U:\W\rightarrow \W$ is said to be
$\nu$-a.s. right-invertible if and only if there exists a measurable map
$V:\W\rightarrow \W$ such that $V\nu\ll\nu$ and  $U\circ V=I_\W$ $\nu$-a.s.
\end{definition}

\begin{proposition}
\label{lrinv}
Set $U,V:\W\rightarrow\W$ measurable maps such that $V\circ U=I_\W$
$\nu$-a.s. and $V\nu\ll\nu$ Then  $U\circ V=I_\W$ $U\nu$-a.s., so if $U\nu\sim \nu$, we
also have $U\circ V=I_\W$ $\nu$-a.s. In that case, we will say that
$U$ is $\nu$-a.s. invertible and we also have $V\nu\sim\nu$.
\end{proposition}

\nproof There exists $A\subset W$ such that $\nu(A)=1$ and for every
$w\in A$, $V\circ U(w)=w$. Consider such a set A, we have
\beaa \EE_{U\nu}\left[1_{U\circ
    V(w)=w}\right]&=&\EE_\nu\left[1_{U\circ V\circ U(w)=U(w)}\right]\\
&=&\EE_\nu\left[1_{U\circ V\circ U(w)=U(w)}1_{w\in
      A}\right]+\EE_\nu\left[1_{U\circ V\circ U(w)=U(w)}1_{w\notin
        A}\right]\\
&=&\EE_\nu\left[1_{ U(w)=U(w)}1_{w\in
      A}\right]\\
&=&1\eeaa
\no Now assume that U is $\nu$-a.s. invertible, set $A\in\F_1$ such that $V\nu(A)=0$. We have $1_A\circ V=0$ $\nu$-a.s. and since $U\nu\sim\nu$, we have $1_A\circ V=0$ $U\nu$-a.s. Finally,
\beaa \nu(A)=\EE_\nu\left[1_A\right]=\EE_\nu\left[1_A\circ V\circ U\right]\eeaa
\no which concludes the proof.
\nqed

\begin{theorem}
\label{inv1}
\no Set $u\in\D$ and assume there exists
a measurable map $V:\W\rightarrow \W$ such that $V\circ W^u=I_\W$
$\nu$-a.s. Denote $v=-u\circ V$. Then $W^u\circ V=I_\W$ $\nu$-a.s.,
$v\in \D$ and $V=W^v$ $\nu$-a.s.\newline
\no Moreover, we have
\beaa \frac{dW^u\nu}{d\nu}&=&\rho(-\delta_\beta v)\\
\frac{dW^v\nu}{d\nu}&=&\rho(-\delta_\beta u)\eeaa
\end{theorem}

\nproof Set $f\in C_b(\W)$, we have
$$\EE_\nu\left[f\circ V\right]=\EE_\nu\left[f\circ V\circ
  W^u\rho(-\delta_\beta u)\right]=\EE_\nu\left[f\rho(-\delta_\beta
  u)\right]$$
\no So $V\nu\sim\nu$ and
$$\frac{dV\nu}{d\nu}=\rho(-\delta_\beta u)$$
\no Since $W^u\nu\sim\nu$, the first assertion comes from proposition \ref{lrinv}.

\no Clearly $v\in L^0(\nu,H)$. Since $u\in \D$,
there exists $n\in\NN$ such that $d\nu\times dt$-a.s.
$$\left|\dot{u}(s,w)\right|\leq n$$
\no Since $V\nu\ll\nu$, we have $d\nu\times dt$-a.s.
$$\left|\dot{v}(s,w)\right|\leq n$$
\no Finally, let us prove that $\dot{v}$ is adapted. We have $\nu$-a.s.
$$\dot{v}\circ W^u=-\dot{u}\circ
V\circ W^u=-\dot{u}$$
\no hence $v\circ W^u$ is adapted. Set $A\in L^2(d\nu\times dt)$ an adapted process. We have:
\beaa\EE_\nu\left[\rho(-\delta_\beta u)\int_0^1\dot{v}(s)\circ W^u
  A(s)\circ W^uds\right]&=&\EE_\nu\left[\int_0^1\dot{v}(s)
  A(s)ds\right]\\
&=&\EE_\nu\left[\int_0^1\EE_\nu\left[\dot{v}(s)|\F_s\right]
  A(s)ds\right]\\
&=&\EE_\nu\left[\rho(-\delta_\beta u)\int_0^1\EE_\nu\left[\dot{v}(s)|\F_s\right]\circ W^u
  A(s)\circ W^uds\right]\eeaa
\no So $\EE_\nu\left[\dot{v}(s)|\F_s\right]\circ W^u=\dot{v}(s)\circ
W^u$ $ds\times d\nu$-a.s. which implies
$\EE_\nu\left[\dot{v}(s)|\F_s\right]=\dot{v}(s)$
$ds\times d\nu$-a.s. since $W^u\nu\sim \nu$.\newline
\no We have
$$W^v\circ W^u=W^{u+v\circ W^u}=W^0=I_\W\;\;\;\nu-a.s.$$
\no and
$$V=W^v\circ W^u\circ V=W^v\;\;\;\nu-a.s.$$
Finally, set $f\in C_b(\W)$,
\beaa \EE_\nu\left[f\circ W^u\right]&=&\EE_\nu\left[f\circ W^u\circ
  W^v\rho(-\delta_\beta v)\right]\\
&=&\EE_\nu\left[f\rho(-\delta_\beta v)\right]\\
\EE_\nu\left[f\circ W^v\right]&=&\EE_\nu\left[f\circ W^v\circ
  W^u\rho(-\delta_\beta u)\right]\\
&=&\EE_\nu\left[f\rho(-\delta_\beta u)\right]\eeaa
\no which gives the final assertion.\nqed

\begin{remark} Set $u\in \D$, $W^u$ is $\nu$-a.s. invertible if and only if
  it is $\nu$-a.s. left-invertible.
\end{remark}

\subsection{\bf{Entropic characterisation of the invertibility of $W^u$}\label{si}}

\no In this section, we prove that the process $W^u$ is left invertible if
and only if the kinetic energy of the perturbation u is equal to
the relative entropy of $W^u\nu$.

\begin{proposition}
\label{i1}
\no Set $u\in  \D$. We have:
$$H(W^u\nu|\nu)\leq\frac{1}{2}\EE_{\nu}\left[|u|^2_H\right]$$
\end{proposition}

\nproof Set $g\in C_b(\W)$ and denote $L=\frac{dW^u\nu}{d\nu}$, we
have:
\beaa \EE_\nu\left[g\circ W^u\right]&=&\EE_{\nu}\left[gL\right]\\
&=&\EE_{\nu}\left[g\circ W^u L\circ W^u\rho(-\delta_\beta u)\right]\eeaa
\no Hence $L\circ W^u\EE_\nu\left[\left.\rho(-\delta_\beta u)\right|\F_1^{W^u}\right]= 1$
$\nu$-a.s. and

\beaa H(W^u\nu|\nu)&=&\EE_{\nu}\left[L\log L\right]\\
&=&\EE_{W^u\nu}\left[\log L\right]\\
&=&\EE_\nu\left[\log L\circ W^u\right]\\
&=&-\EE_\nu\left[\log\EE_\nu\left[\left.\rho(-\delta_\beta u)\right|\F_1^{W^u}\right]\right]\\
&\leq&-\EE_\nu\left[\log\rho(-\delta_\beta u)\right]\\
&\leq&\frac{1}{2}\EE_\nu\left[|u|_H^2\right]\eeaa
\nqed
\newline

\no The proof gives the following additional result

\no\begin{corollary}
For $u\in\D$, we have
$$L\circ W^u\EE_\nu\left[\left.\rho(-\delta_\beta u)\right|\F_1^{W^u}\right]=1$$
\end{corollary}
\no Now comes the criteria:

\begin{theorem}
\label{inv}
Set $u\in\D$, then $W^u$ is $\nu$-a.s. invertible
if and only if:
$$H(W^u\nu|\nu)=\frac{1}{2}\EE_{\nu}\left[|u|_H^2\right]$$
\end{theorem}

\nproof Assume that the inequality hold. We still denote
$L=\frac{dW^u\nu}{d\nu}$ and as in last proof we have $\nu$-a.s.
$$L\circ W^u\EE_\nu\left[\left.\rho(-\delta_\beta u)\right|\F_1^{W^u}\right]= 1$$

\no Using Jensen inequality we have $\nu$-a.s.

\beaa 0&=& \log L\circ
 W^u+\log\EE_\nu\left[\left.\rho(-\delta_{\beta}
    u)\right|\F^{ W^u}_1\right]\\
&\geq& \log L\circ
 W^u+\EE_\nu\left[\left. \log\rho(-\delta_{\beta}
    u)\right|\F^{ W^u}_1\right]\eeaa
\no and
\beaa 0&\geq&\EE_\nu\left[\log L\circ W^u\right]+\EE_\nu\left[\log\rho(-\delta_{\beta}
    u)\right]\\
&\geq&H( W^u\nu|\nu)-\frac{1}{2}\EE_\nu\left[|u|_H^2\right]\\
&=&0\eeaa
\no So
\beaa 0&=&\log L\circ W^u+\log\EE_\nu\left[\left.\rho(-\delta_{\beta} u)\right|\F_1^{ W^u}\right]\\
&=&\log L\circ  W^u+\EE_\nu\left[\log\left.\rho(-\delta_{\beta} u)\right|\F_1^{ W^u}\right]\eeaa
\no and
$$\log\EE_\nu\left[\left.\rho(-\delta_{\beta} u)\right|\F_1^{ W^u}\right]=\EE_\nu\left[\log\left.\rho(-\delta_{\beta} u)\right|\F_1^{ W^u}\right]$$
\no The strict concavity of the function $\log$ gives
$$\EE_\nu\left[\left.\rho(-\delta_{\beta} u)\right|\F_1^{ W^u}\right]=\rho(-\delta_{\beta} u)$$
\no Finally we have
\begin{equation}\label{2eg}L\circ W^u\rho(-\delta_\beta u)= 1\end{equation}
\no Since $\beta$ is a $\nu$-Brownian motion, there exists $v\in
L^0_a(\nu,H)$ such that
$L=\rho(-\delta_\beta v)$.
\newline
\no We apply the logarithm to \ref{2eg} to get:
\beaa 0&=&\delta_\beta v\circ W^u+\frac{1}{2}|v\circ W^u|_H+\delta_\beta
u+ \frac{1}{2}|u|_H\eeaa
\no We have:
\beaa \delta_\beta v\circ W^u &=&\int_0^1\dot{v}(s)\circ W^ud\beta(s) +\langle v\circ W^u,u\rangle_H\eeaa

\no so finally we have:
\begin{equation}\label{2test} 0= \delta_\beta(v\circ W^u+u)+\frac{1}{2}|v\circ W^u+u|_H^2\end{equation}
According to Girsanov theorem $\beta+v$ is a $W^u\nu$-Brownian motion,
so:
\beaa \EE_{\nu}\left[L\log L\right] &=& \EE_{W^u\nu}\left[\log
    L\right]\\
&=&\EE_{W^u\nu}\left[-\int_0^1\dot{v}(s)d\beta(s)-\frac{1}{2}\int_0^1|\dot{v}(s)|^2ds\right]\\
&=&\frac{1}{2}\EE_{W^u\nu}\left[\int_0^1|\dot{v}(s)|^2ds\right]\\
&=&\frac{1}{2}\EE_\nu\left[|v\circ W^u|_H^2\right]\eeaa

\no So $v\circ W^u\in L^2_a(\nu,H)$ and we can take the expectation
with respect to $\nu$ in \ref{2test} to obtain $u+v\circ W^u=0$
$\nu$-a.s. which
implies $v\in \D$. So $\nu$-a.s.
$$W^v\circ W^u=W^{u+v\circ W^u}=W^0=I_\W$$
\no Conversely, assume that $W^u$ admits an inverse V and set
$v=-u\circ V$. According to theorem \ref{inv1}, $v\in D$ and $W^v=V$
$\nu$-a.s.
\no Once again, denote $L=\frac{dW^u\nu}{d\nu}$. We know that
$L=\rho(-\delta_\beta v)$. Observe that
\beaa \log L\circ W^u &=&
\left(-\int_0^1\dot{v}(s)d\beta(s)-\frac{1}{2}\int_0^1\dot{v}(s)^2ds\right)\circ
W^u\\
&=&-\log\rho(-\delta_\beta u)\eeaa
\no So finally
\beaa H(W^u\nu|\nu)&=&\EE_\nu[L\log L] =\EE_\nu[\log L\circ W^u]\\
&=&\EE_\nu\left[-\log\rho(-\delta_\beta u)\right]\\
&=&\frac{1}{2}\EE_\nu[|u|_H^2]\eeaa
\nqed
\newline
\no The proof gives the following additional result:
\begin{corollary} Set
  $u\in\D$
such that  $W^u$ is $\nu$-a.s. left-invertible, we have
$$L\circ W^u\rho(-\delta_\beta
    u)=1$$
\end{corollary}

\no The following corollary is an immediate consequence of theorems \ref{inv1} and \ref{inv}.

\begin{corollary}
Set $u\in\D$, we have
$$H(W^u\nu|\nu)=\frac{1}{2}\EE_\nu\left[|u|_H^2\right]$$
\no if and only if there exists $v\in\D$ such that
$$W^v\circ W^u=W^u\circ W^v=I_\W\;\nu-a.s.$$
\end{corollary}

\begin{definition}
\no We define $\D^{i}$ as
$$\D^i=\left\{u\in D,W^u\;is\;\nu-a.s.\;invertible\right\}$$
\end{definition}

\section{\bf{Variational problem}}
\subsection{\bf{Approximation of absolutely continuous measures}}

\begin{theorem}
\label{repr}
If $\theta\sim \nu$ is such that there exists $r>1$ such that
$$\frac{d\theta}{d\nu}\log\frac{d\theta}{d\nu}\in L^1(\nu)$$
\no and
$$\log\frac{d\theta}{d\nu}\in L^r(\nu)$$
\no there exists $(u_n)\in \left(\D^i\right)^\NN$ such that,
\beaa&&\frac{dW^{u_n}\nu}{d\nu}\log \frac{dW^{u_n}\nu}{d\nu}\rightarrow\frac{d\theta}{d\nu}\log
\frac{d\theta}{d\nu}\;\; in\;\;L^1(\nu)\\
&&\frac{dW^{u_n}\nu}{d\nu}\log \frac{d\theta}{d\nu}\rightarrow\frac{d\theta}{d\nu}\log
\frac{d\theta}{d\nu}\;\; in\;\;L^1(\nu)\eeaa
\end{theorem}

\nproof Denote
\beaa L&=&\frac{d\theta}{d\nu}\eeaa
\no Eventually sequentializing afterward, we have to prove that for any
$\epsilon>0$, there exists $u\in\D^i$ such that
\beaa \EE_\nu\left[\left|\frac{dW^u\nu}{d\nu}\log \frac{dW^u\nu}{d\nu}-L_1\log
    L_1\right|\right]&\leq&\epsilon\\
\EE_\nu\left[\left|\frac{dW^u\nu}{d\nu}\log L_1-L_1\log
    L_1\right|\right]&\leq&\epsilon\eeaa

\no The proof is divided in six steps.\newline
\no Step 1 : We approximate L with a density that is both lower-bounded and upper bounded.\newline
\no Denote
\beaa\phi_n&=&\min(L,n)\\
L_n&=&\frac{\phi_n}{\EE_\nu\left[\phi_n\right]}\eeaa
\no The monotone convergence theorem ensures that
$\EE_\nu\left[\phi_n\right]\rightarrow 1$ so for any $\alpha\in (0,1)$,
there exists some $n_\alpha\in\NN$ such that for any $n\geq n_\alpha$,
$$\EE_\nu\left[\phi_n\right]\geq\alpha$$
\no $(L_n\log L_n)$ converges $\nu$-a.s. to $L\log L$ and if $n\geq n_\alpha
$ and
\beaa
\left|L_n\log
  L_n\right|&=&\left|\frac{\phi_n}{\EE_\nu\left[\phi_n\right]}\log
    \frac{\phi_n}{\EE_\nu\left[\phi_n\right]}\right|1_{\frac{\phi_n}{\EE_\nu\left[\phi_n\right]}\leq 1}+
\left|\frac{\phi_n}{\EE_\nu\left[\phi_n\right]}\log
    \frac{\phi_n}{\EE_\nu\left[\phi_n\right]}\right|1_{\frac{\phi_n}{\EE_\nu\left[\phi_n\right]}>
    1}\\
&\leq& e^{-1}1_{\frac{\phi_n}{\EE_\nu\left[\phi_n\right]}\leq 1}+
\left|\frac{L}{\alpha}\log
   \frac{L}{\alpha}\right|1_{\frac{\phi_n}{\EE_\nu\left[\phi_n\right]}>
    1}\\
&\leq& e^{-1}+\left|\frac{L}{\alpha}\log
    \frac{L}{\alpha}\right|\eeaa

\no So the Lebesgue theorem ensures that $(L_n\log L_n)$ converge
toward $L\log L$ in $L^1(\nu)$.
\no Similarly, $(L_n\log L)$ converges $\nu$-a.s. to $L\log L$ and
if $n\leq n_\alpha$,
\beaa\left|L_n\log L\right|\leq \left|\frac{L}{\alpha}\log L\right|\eeaa
\no and the Lebesgue theorem ensures that $(L_n\log L)$ converges to
$L_n\log L$ in $L^1(\nu)$, so there exists $n_0\in\NN$ such that
\beaa \EE_\nu\left[\left|L_{n_0}\log L_{n_0}-L\log
    L\right|\right]&\leq&\epsilon\\
\EE_\nu\left[\left|L_{n_0}\log
    L-L\log L\right|\right]&\leq&\epsilon\eeaa
\no $\left(\frac{L_{n_0}+a}{1+a}\log \frac{L_{n_0}+a}{1+a}\right)$ converges
$\nu$-a.s. to $L_{n_0}\log L_{n_0}$ when a converges to 0. Set $a\in [0,1]$, we have
\beaa \left|\frac{L_{n_0}+a}{1+a}\log \frac{L_{n_0}+a}{1+a}\right|&=&\left|\frac{L_{n_0}+a}{1+a}\log
\frac{L_{n_0}+a}{1+a}\right|1_{L_{n_0}\leq 1}+\left|\frac{L_{n_0}+a}{1+a}\log
\frac{L_{n_0}+a}{1+a}\right|1_{L_{n_0}>1}\\
&\leq& e^{-1}1_{L_{n_0}\leq 1}+\left|L_{n_0}\log
  L_{n_0}\right|1_{L_{n_0}>1}\\
&\leq& e^{-1}+\left|L_{n_0}\log
  L_{n_0}\right|\eeaa
So the Lebesgue theorem ensures that $\left(\frac{L_{n_0}+a}{1+a}\log \frac{L_{n_0}+a}{1+a}\right)$ converges
to $L_{n_0}\log L_{n_0}$ in $L^1(\nu)$. Similarly,
$\left(\frac{L_{n_0}+a}{1+a}\log L\right)$ converges $\nu$-a.s. to
$L_{n_0}\log L$ and
\beaa \left|\frac{L_{n_0}+a}{1+a}\log L\right|&\leq&\left|(L_{n_0}+1)\log L\right|\eeaa
\no and the Lebesgue theorem ensures that $\left(\frac{L_{n_0}+a}{1+a}\log L\right)$ converges to
$L_{n_0}\log L$ in $L^1(\nu)$ and
there exists $a\in [0,1]$ such that
\beaa \EE_\nu\left[\left|\frac{L_{n_0}+a}{1+a}\log
  \frac{L_{n_0}+a}{1+a}-L_{n_0}\log L_{n_0}\right|\right]\leq\epsilon\\
\EE_\nu\left[\left|\frac{L_{n_0}+a}{1+a}\log
  L-L_{n_0}\log L\right|\right]\leq\epsilon\eeaa

\no $\frac{L_{n_0}+a}{1+a}$ is both lower-bounded and upper-bounded in
$L^\infty(\nu)$,
denote these bounds respectively d and D.\newline
\no Also denote
$$M(t)=\EE_\nu\left[\left.\frac{L_{n_0}+a}{1+a}\right|\F_t\right]$$
\no We write
$$M(t)=\exp\left(\int_0^t\dot{\alpha}(s)d\beta(s)-\frac{1}{2}\int_0^t\left|\dot{\alpha}(s)\right|^2ds\right)$$
\no with $\alpha\in
L^0_a(\nu,H)$.\newline
\no Step 2 : We prove that $\alpha\in L^2(\nu,H)$.\newline
\no Set
$$T_n=\inf\left\{t\in [0,1], \int_0^t\left|\dot{\alpha}(s)\right|^2ds>n\right\}$$
\no $(T_n)$ is a sequence of stopping times which increases
stationarily toward 1. We have, using
$M=1+\int_0^.\dot{\alpha}(s)M(s)d\beta(s)$
\beaa \EE_\nu\left[\left(M(t\wedge
      T_n)-1\right)^2\right]&=&\EE_\nu\left[\int_0^{t\wedge T_n}\left|\dot{\alpha}(s)\right|^2M(s)^2ds\right]\\
&\geq& d^2\EE_\nu\left[\int_0^{t\wedge T_n}\left|\dot{\alpha}(s)\right|^2ds\right]\eeaa
\no so
$$\EE_\nu\left[\int_0^{t\wedge T_n}\left|\dot{\alpha}(s)\right|^2ds\right]\leq \frac{1}{d^2}\EE_\nu\left[\left(M(t\wedge
      T_n)-1\right)^2\right]\leq \frac{2\left(D^2+1\right)}{d^2}$$
\no hence passing to the limit
$$\EE_\nu\left[\int_0^{1}\left|\dot{\alpha}(s)\right|^2ds\right]\leq\infty$$
\no Step 3 : we approximate $\alpha$ with an element of $L^\infty(\nu,H)$.\newline
\no Define
$$\alpha_n(t,w)\in\RR\times W\mapsto\int_0^t
\dot{\alpha}(s,w)1_{[0,T_n]}(s,w)ds$$

\no and
\beaa
M^n(t)&=&\exp\left(\int_0^t\dot{\alpha^n}(s)d\beta(s)-\frac{1}{2}\int_0^t\left|\dot{\alpha^n}(s)\right|^2ds\right)\eeaa
\no and clearly $(M^n(1)\log M^n(1))$ converges $\nu$-a.s. to $M(1)\log M(1)$, $(M^n(1)\log
L)$ converges $\nu$-a.s. to $M(1)\log L$ and
$M^n(1)=\EE_\nu\left[\left. M(1)\right|\F_{T_n}\right]$, so $\nu$-a.s.
\beaa \left|M^n(1)\log M^n(1)\right|&\leq &\max \left(e^{-1},\left|D\log
  D\right|\right)\\
\left|M^n(1)\log L\right|&\leq&\left|D\log L\right|\eeaa
\no so the Lebesgue theorem ensures that $(M^n(1)\log M^n(1))$ converges
to $M(1)\log M(1)$ in $L^1(\nu)$ and $(M_1^n\log
L)$ converges to $M_1\log L$ in $L^1(\nu)$ and there exists
$n\in\NN$ such that
\beaa \left|M^{n}(1)\log M^{n}(1)-M(1)\log M(1)\right|&\leq&\epsilon\\
\left|M^{n}(1)\log L-M(1)\log L\right|&\leq&\epsilon\eeaa
\no Step 4 : We approximate $\alpha^n$ with an element of $\D$.\newline
\no Define $$\xi^{n,m}:(t,w)\in [0,1]\times\W\mapsto\int_0^t\max \left(\min
  \left(\dot{\alpha}^n(s,w),m\right),-m\right)ds$$
\no and
$$
M^{n,m}(t)=\exp\left(\int_0^t\dot{\xi^{n,m}}(s)d\beta(s)-\frac{1}{2}\int_0^t\left|\dot{\xi^{n,m}}(s)\right|^2ds\right)$$
\no $\left(M^{n,m}(1)\log M^{n,m}(1)\right)$ and $\left(M^{n,m}(1)\log
  L\right)$ converges respectively to $M^n(1)\log M^n(1)$ and $M^n(1)\log
L$ in probability. To prove that $\left(M^{n,m}(1)\log
  M^{n,m}(1)\right)$ is uniformly integrable, it is sufficient to prove
it is bounded in any $L^p(\nu)$, set $p>1$
\beaa
\EE_\nu\left[\left|M^{n,m}(1)\right|^{p}\right]&=&\EE_\nu\left[\exp\left(p\int_0^1\dot{\xi^{n,m}}(s)d\beta(s)-\frac{p}{2}\int_0^1\left|\dot{\xi^{n,m}}(s)\right|^2ds\right)\right]\\
&=&\EE_\nu\left[\exp\left(p\int_0^1\dot{\xi^{n,m}}(s)d\beta(s)-\frac{p^2}{2}\int_0^1\left|\dot{\xi^{n,m}}(s)\right|^2ds\right)\exp\left(\frac{p^2-p}{2}\int_0^1\left|\dot{\xi^{n,m}}(s)\right|^2ds\right)\right]\\
&\leq&\EE_\nu\left[\exp\left(\int_0^1p\dot{\xi^{n,m}}(s)d\beta(s)-\frac{1}{2}\int_0^1\left|p\dot{\xi^{n,m}}(s)\right|^2ds\right)\exp\left(\frac{p^2-p}{2}n\right)\right]\\
&\leq&\exp\left(\frac{p^2-p}{2}n\right)\eeaa
\no so $\left(M^{n,m}(1)\log M^{n,m}(1)\right)$ converges to $M^n(1)\log
M^n(1)$ in $L^1(\nu)$. Furthermore, set p such that
$p^{-1}+r^{-1}=1$
\beaa \EE_\nu\left[\left|M^{n,m}(1)\log L-M^n(1)\log
    L\right|\right]&\leq&\left|M^{n,m}(1)-M^n(1)\right|_{L^p(\nu)}\left|\log
  L\right|_{L^r(\nu)}\\
&\rightarrow &0\eeaa
\no and there exists some $m>0$ such that
\beaa \EE_\nu\left[\left|M^{n,m}(1)\log M^{n,m}(1)-M^n(1)\log
M^n(1)\right|\right]&\leq&\epsilon\\
\EE_\nu\left[\left|M^{n,m}(1)\log L-M^n(1)\log L\right|\right]&\leq&\epsilon\eeaa
\no Step 5 : we approximate $\xi^{n,m}$ with a retarded shift.\newline
\no For $\eta>0$ set
$$\dot{\gamma^\eta}(t,w)\in[0,1]\times W\mapsto
\dot{\xi^{n,m}}(t-\eta)(w)1_{t>\eta}$$
$$N^{\eta}(t)=\exp\left(\int_0^1\dot{\gamma^\eta}(s)d\beta(s)-\frac{1}{2}\int_0^1\left|\dot{\gamma^\eta}(s)\right|^2ds\right)$$
\no We have for any $\eta>0$, $d\nu\times ds$-a.s.,
$$\left|\dot{\gamma^\eta}(s)\right|\leq m$$
\no i.e. $\gamma^\eta\in\D$.\newline
\no Similarly as before $\left(N^\eta(1)\log N^\eta(1)\right)$
converges in probability to $
M^{n,m}(1)\log M^{n,m}(1)$ and $\left(N^\eta(1)\right)$ is bounded in every
$L^p(\nu)$ hence $\left(N^\eta(1)\log N^\eta(1)\right)$ is uniformly
integrable and converges in $L^1(\nu)$ to $
M^{n,m}(1)\log M^{n,m}(1)$\newline
\no Furthermore, using Holder inequality, we have
$$\EE_\nu\left[\left|N^\eta(1)\log L-M^{n,m}(1)\log
    L\right|\right]\leq\left|N^\eta(1)-M^{n,m}(1)\right|_{L^{p}(\nu)}\left|\log
  L\right|_{L^r(\nu)}$$
\no where $\frac{1}{p}+\frac{1}{r}=1$.\newline
\no Consequently there exists $\eta>0$ such that
\beaa \EE_\nu\left[\left|N^\eta(1)\log N^\eta(1)-M^{n,m}(1)\log
    M^{n,m}(1)\right|\right]&\leq&\epsilon\\
 \EE_\nu\left[\left|N^\eta(1)\log L-M^{n,m}(1)\log
    L\right|\right]&\leq&\epsilon\eeaa
\no using triangular inequality, we have
\beaa \EE_\nu\left[\left|L\log L-N^\eta(1)\log
  N^\eta(1)\right|\right]&\leq&\EE_\nu\left[\left|L\log L-L_{n_0}\log
L_{n_0}\right|\right]\\
&&+\EE_\nu\left[\left|L_{n_0}\log
L_{n_0}-\frac{L_{n_0}+a}{1+a}\log\frac{L_{n_0}+a}{1+a}\right|\right]\\
&&+\EE_\nu\left[\left|\frac{L_{n_0}+a}{1+a}\log\frac{L_{n_0}+a}{1+a}-M^n(1)\log
M^n(1)\right|\right]\\
&&+\EE_\nu\left[\left|M^n(1)\log M^n(1)-M^{n,m}(1)\log
M^{n,m}(1)\right|\right]\\
&&+\EE_\nu\left[\left|M^{n,m}(1)\log M^{n,m}(1)-N^\eta(1)\log
N^\eta(1)\right|\right]\\
&\leq&5\epsilon\\
\EE_\nu\left[\left|L\log L-N^\eta(1)\log
 L\right|\right]&\leq&\EE_\nu\left[\left|L\log L-L_{n_0}\log
L\right|\right]\\
&&+\EE_\nu\left[\left|L_{n_0}\log
L-\frac{L_{n_0}+a}{1+a}\log L\right|\right]\\
&&+\EE_\nu\left[\left|\frac{L_{n_0}+a}{1+a}\log L-M^n(1)\log
L\right|\right]\\
&&+\EE_\nu\left[\left|M^n(1)\log L-M^{n,m}(1)\log
L\right|\right]\\
&&+\EE_\nu\left[\left|M^{n,m}(1)\log L-N^\eta(1)\log
L\right|\right]\\
&\leq&5\epsilon\eeaa

\no Step 6 : We prove that $ W^{-\gamma^{\eta}}$ is
$\nu$-a.s. left-invertible and is the solution to our problem.\newline
\no Set $A\subset\W$ such that $\nu(A)=1$ and for every $w\in A$, $\beta\circ W^{-\gamma^\eta}(w)=\beta(w)-\gamma^\eta(w)$ and set $w_1,w_2\in A$ such that
$ W^{-\gamma^\eta}(w_1)= W^{-\gamma^\eta}(w_2)$. We have
\beaa \beta\circ W^{-\gamma^\eta}(w_1)&=&\beta\circ W^{-\gamma^\eta}(w_2)\\
\beta(w_1)-\int_0^.\dot{\gamma^\eta}(s,w_1)ds&=&\beta(w_2)-\int_0^.\dot{\gamma^\eta}(s,w_2)ds\eeaa

\no For any $s\in[0,\eta]$, $\beta(s,w_1)=\beta(s,w_2)$,
$\gamma^\eta$ being adapted to filtration $(\F^\beta_{s-\eta})$,
it implies that for $s\in[0,2\eta]$
$$\int_0^s\dot{\gamma^\eta}(r,w_1)ds=\int_0^s\dot{\gamma^\eta}(r,w_2)ds$$
\no and
$$\beta(s,w_1)=\beta(s,w_2)$$
\no An easy iteration shows that $\beta(w_1)=\beta(w_2)$.\newline

\no Since $\beta$ and W have the same filtrations, and $\beta$ is $\nu$-a.s. path-continuous, we can write $W(t)=\phi_t(\beta(s),s\in[0,t]\cap\QQ)$ $\nu$-a.s. for every $t\in [0,1]$, with $\phi$ a measurable function from $\RR^QQ$ to $\RR$, see \cite{nev}. Consequently, we can write $\left(W(t),t\in[0,1]\cap\QQ\right)=\phi\left(\beta(t), t\in[0,1]\cap\QQ\right)$ $\nu$-a.s., with $\phi$ a measurable function from $\RR^\QQ$ to $\RR^\QQ$. Denote
\beaa A'=A\cap\left\{w\in\W, \left(W(t,w),t\in[0,1]\cap\QQ\right)=\phi\left(\beta(t,w), t\in[0,1]\cap\QQ\right)\right\}\eeaa
\no $\nu(A')=1$. Set $w_1,w_2\in A'$ such that $ W^{-\gamma^\eta}(w_1)= W^{-\gamma^\eta}(w_2)$. We have $\beta(w_1)=\beta(w_2)$ so
\beaa\left(W(t,w_1),t\in[0,1]\cap\QQ\right)&=&\left(W(t,w_2),t\in[0,1]\cap\QQ\right)\\
\left(w_1(t),t\in[0,1]\cap\QQ\right)&=&\left(w_2(t),t\in[0,1]\cap\QQ\right)\eeaa

\no $w_1$ and $w_2$ are continuous and coincide on $[0,1]\cap\QQ$ so they are equal.

\no $ W^{-\gamma^{\eta}}$ is $\nu$-a.s. injective and so
$\nu$-a.s. left-invertible, its inverse is of
the form $ W^{v^{\eta}}$, with $v^\eta\in\D$ and we have
$$\frac{d  W^{v^\eta}\nu}{d\nu}=L^{\eta,n}_1$$
\no So $ W^{v^\eta}\nu\sim \nu$ and
$$ W^{v^\eta}\circ
 W^{-\gamma^\eta}=  W^{-\gamma^\eta}\circ  W^{v^\eta}\;\;
\nu-a.s.$$
\nqed

\subsection{\bf{Main theorem}}
\no As stated in the beginning, we aim to provide a variational
representation of $-\log\EE_\nu\left[e^{-f}\right]$ . This first
result is from \cite{art}:

\begin{theorem}
Set $f:\W\rightarrow \RR$ a measurable function verifying
$$\EE_\nu\left[|f|(1+e^{-f})\right]<\infty$$
\no Denote $\mathcal{P}$ the set of probability measures on $\left(\W,\F\right)$ which are absolutely continuous with respect to $\nu$, then
$$-\log\EE_\nu\left[e^{-f}\right]=\inf_{\theta\in
  \mathcal{P}}\left(\EE_\theta[f]+H(\theta|\nu)\right)$$
\no and the unique supremum is attained at the measure
$$d\theta_0=\frac{e^{-f}}{\EE_\nu\left[e^{-f}\right]}d\nu$$
\end{theorem}

\begin{proposition}
\label{ineqvar}
\no Set $f:\W\rightarrow\RR$ a measurable function
verifying $\EE_\nu\left[|f|(1+e^{-f})\right]<\infty$, then
$$-\log\EE_\nu\left[e^{-f}\right]\leq\inf_{u\in\D}\EE_\nu\left[f\circ W^u+\frac{1}{2}|u|_H^2\right]$$
\end{proposition}

\nproof Set $u\in\D$
\beaa
-\log\EE_\nu\left[e^{-f}\right]&\leq&\EE_{W^u\nu}\left[f\right]+H(W^u\nu|\nu)\\
&=&\EE_\nu\left[f\circ W^u+\frac{1}{2}|u|_H^2\right]\eeaa\nqed

\no Here is the main result.

\begin{theorem}
\label{tvar}
Set $p>1$ and $f\in L^p(\nu)$ such that
$\EE_\nu\left[(|f|+1)e^{-f}\right]<\infty$, then we have
$$-\log\EE_\nu\left[e^{-f}\right]=\inf_{u\in\D^i}\EE_\nu\left[f\circ W^u+\frac{1}{2}|u|_H^2\right]$$
\end{theorem}

\nproof Using proposition \ref{ineqvar}, we have easily
$$-\log\EE_\nu\left[e^{-f}\right]\leq\inf_{u\in\D^i}\EE_\nu\left[f\circ W^u+\frac{1}{2}|u|_H^2\right]$$
\no Let $\theta_0$ be the measure on $\W$ defined by
$$d\theta_0=\frac{e^{-f}}{\EE_\nu\left[e^{-f}\right]}d\nu$$
\no According to theorem \ref{repr}, there exists $(u_n)\in
\D^i$ such that for every $n\in\NN$
\beaa \frac{dW^{u_n}\nu}{d\nu}\log \frac{dW^{u_n}\nu}{d\nu} \rightarrow
\frac{d\theta_0}{d\nu}\log \frac{d\theta_0}{d\nu}\\
\frac{dW^{u_n}\nu}{d\nu}\log \frac{d\theta_0}{d\nu} \rightarrow
\frac{d\theta_0}{d\nu}\log \frac{d\theta_0}{d\nu}\eeaa
\no in $L^1(\nu)$.\newline
\no Since $W^{u_n}$
is $\nu$-a.s. invertible, we have
\beaa \EE_\nu\left[f\circ W^{u_n}+\frac{1}{2}|u_n|_H^2\right]=
\EE_\nu\left[f \frac{dW^{u_n}\nu}{d\nu}\right]+\EE_\nu\left[\frac{dW^{u_n}\nu}{d\nu}\log
  \frac{dW^{u_n}\nu}{d\nu}\right]\eeaa
\no When n goes to infinity, we have
$$\EE_\nu\left[\frac{dW^{u_n}\nu}{d\nu}\log
  \frac{dW^{u_n}\nu}{d\nu}\right]\rightarrow
\EE_\nu\left[\frac{d\theta_0}{d\nu}\log \frac{d\theta_0}{d\nu}\right]$$
\no and since $f=-\log \frac{d\theta_0}{d\nu}-\log\EE_\nu\left[e^{-f}\right]$,
\beaa\EE_\nu\left[f\frac{dW^{u_n}\nu}{d\nu}\right]&\rightarrow& \EE_\nu\left[f\frac{d\theta_0}{d\nu}\right]\eeaa
\no So finally, when n goes to infinity,
\beaa \EE_\nu\left[f\circ
  W^{u_n}+\frac{1}{2}|u_n|_H^2\right]&\rightarrow&\EE_{\theta_0}\left[f\right]+H(\theta_0|\nu)\\
&=&-\log\EE_\nu\left[e^{-f}\right]\eeaa
\no which conclude the proof.\nqed

\subsection{\bf{Retrieving the Pr\'ekopa-Leindler theorem}}

\begin{definition}
We denote
$$H_b=\left\{h\in H, \dot{h}\;is\;dt-a.s.\;bounded\right\}$$
\end{definition}

\begin{remark}
Observe that $H_b\subset\D$ and that if $u\in\D$, $u(w)\in H_b$ $\nu$-a.s.
\end{remark}

\begin{theorem}
\label{pl}
Assume that for any $u\in\D$,
$$W^{u}(w)=W^{u(w)}(w)\;\nu-a.s.$$
Set $a,b,c:\W\rightarrow\RR$ positive and measurable such that for every $h,k\in H_b$ and $ t\in [0,1]$ we have $\nu$-a.s.
$$a\circ W^{ t h+(1- t)k}\exp\left(-\frac{1}{2}\left| t
  h+(1- t) k\right|_H^2\right)\geq\left(b\circ W^{h}\exp\left(-\frac{1}{2}\left|
  h\right|_H^2\right)\right)^{ t}\left(c\circ W^{k}\exp\left(-\frac{1}{2}\left|
  k\right|_H^2\right)\right)^{1- t}$$
\no then for any density d such that $h\in H_b\mapsto -\log d\circ W^h$ is
$\nu$-a.s. concave, if $\theta$ denotes the measure on $\W$ given by
$\frac{d\theta}{d\nu}=d$, we have in $\bar{\RR}$:
$$\EE_\theta\left[a\right]\geq\left(\EE_\theta\left[b\right]\right)^ t\left(\EE_\theta\left[c\right]\right)^{1- t}$$
\end{theorem}

\nproof First observe that eventually replacing a,b,c with $da,db,dc$
we only need to prove the case $d=1$ i.e. $\theta=\nu$\newline

\no With the convention $\log(\infty)=\infty$ and $\log(0)=-\infty$, we denote
$$\tilde{a}=-\log a,\tilde{b}=-\log b,\tilde{c}=-\log c$$
\no We begin with the case where there exists $m,M>0$ such that we
have $\nu$-a.s.
$$m\leq \tilde{a},\tilde{b},\tilde{c}\leq M$$
\no Set $ t\in [0,1]$, for $h,k\in H_b$, we have $\nu$-a.s.
\beaa &&a\circ W^{ t h+(1- t)k}\exp\left(-\frac{1}{2}\left| t
  h+(1- t) k\right|_H^2\right)\\
  &&\geq\left(b\circ W^{h}\exp\left(-\frac{1}{2}\left|
  h\right|_H^2\right)\right)^{ t}\left(c\circ W^{k}\exp\left(-\frac{1}{2}\left|
  k\right|_H^2\right)\right)^{1- t}\eeaa
  \no So for $u_1,u_2\in\D^i$
  \beaa &&a\circ W^{ t u_1+(1- t)u_2}\exp\left(-\frac{1}{2}\left| t
  u_1+(1- t) u_2\right|_H^2\right)\\
  &&\geq\left(b\circ W^{u_1}\exp\left(-\frac{1}{2}\left|
  u_1\right|_H^2\right)\right)^{ t}\left(c\circ W^{u_2}\exp\left(-\frac{1}{2}\left|
  u_2\right|_H^2\right)\right)^{1- t}\eeaa
\no hence applying the logarithm function , changing the sign and
taking the expectation
\beaa &&\EE_{\nu}\left[\tilde{a}\circ W^{ t u_1 +(1- t)u_2}+\frac{1}{2}\left| t
  u_1+(1- t)u_2\right|^2_H\right]\\
  &&\leq  t\EE_{\nu}\left[\tilde{b}\circ W^{u_1}+\frac{1}{2}\left|u_1\right|^2_H\right]+(1- t)\EE_{\nu}\left[\tilde{c}\circ W^{u_2}+\frac{1}{2}\left|u_2\right|^2_H\right]\eeaa
\no So
\beaa \inf_{u\in \D^i}
\EE_{\nu}\left[\tilde{a}\circ W^u+\frac{1}{2}|u|^2_H\right]\leq
   t\EE_{\nu}\left[\tilde{b}\circ
    W^{u_1}+\frac{1}{2}\left|u_1\right|^2_H\right]+(1- t)\EE_{\nu}\left[\tilde{c}\circ
    W^{u_2}+\frac{1}{2}\left|u_2\right|^2_H\right]\eeaa

\no According to theorem \ref{tvar} we have
\beaa -\log\EE_{\nu}\left[e^{-\tilde{a}}\right]&\leq&
   t\EE_{\nu}\left[\tilde{b}\circ
    W^{u_1}+\frac{1}{2}\left|u_1\right|^2_H\right]+(1- t)\EE_{\nu}\left[\tilde{c}\circ
    W^{u_2}+\frac{1}{2}\left|u_2\right|^2_H\right]\eeaa
\no which implies
\beaa  -\log\EE_{\nu}\left[e^{-\tilde{a}}\right]&\leq&
   t\EE_{\nu}\left[\tilde{b}\circ
    W^{u_1}+\frac{1}{2}\left|u_1\right|^2_H\right]+(1- t)  \inf_{v\D^i} \EE_{\nu}\left[\tilde{c}\circ
    W^{v}+\frac{1}{2}\left|v\right|^2_H\right]\\
&=& t\EE_{\nu}\left[\tilde{b}\circ
    W^{u_1}+\frac{1}{2}\left|u_1\right|^2_H\right]-(1- t)\log\EE_{\nu}\left[e^{-\tilde{c}}\right]\eeaa
\no which implies once again
\beaa  -\log\EE_{\nu}\left[e^{-\tilde{a}}\right]&\leq& \inf_{v\in \D^i}\left(
 t\EE_{\nu}\left[\tilde{b}\circ
    W^{v}+\frac{1}{2}\left|v\right|^2_H\right]\right)-(1- t)\log\EE_{\nu}\left[e^{-\tilde{c}}\right]\\
&=&- t\log\EE_{\nu}\left[e^{-\tilde{b}}\right]
-(1- t)\log\EE_{\nu}\left[e^{-\tilde{c}}\right]\eeaa
\no taking the opposite and applying the exponential, we get
$$\EE_{\nu}\left[e^{-\tilde{a}}\right]\geq\left(\EE_{\nu}\left[e^{-\tilde{b}}\right]\right)^ t\left(\EE_{\nu}\left[e^{-\tilde{c}}\right]\right)^{1- t}$$

\no For the general case, denote for $n\in\NN$ and $m\in\NN^*$
\beaa &&\tilde{a}_n=\tilde{a}\wedge n,\tilde{b}_n=\tilde{b}\wedge n,\tilde{c}_n=\tilde{c}\wedge n\\
&&\tilde{a}_{nm}=\tilde{a}_n+\frac{1}{m},\tilde{b}_{nm}=\tilde{b}_n+\frac{1}{m},\tilde{c}_{nm}=\tilde{c}_n+\frac{1}{m}\eeaa
\no For every $h,k\in H_b$, we have $\nu$-a.s.:
$$\tilde{a}_{nm}\circ W^{ t h +(1- t)k}+\frac{1}{2}\left| t
  h+(1- t)k\right|^2_H\leq  t\tilde{b}_{nm}\circ W^{h}+\frac{1}{2}\left|h\right|^2_H+(1- t)\tilde{c}_{nm}\circ W^{k}+\frac{1}{2}\left|k\right|^2_H$$
\no so the bounded case we treated above ensures that
$$\EE_{\nu}\left[e^{-\tilde{a}_{nm}}\right]\geq\left(\EE_{\nu}\left[e^{-\tilde{b}_{nm}}\right]\right)^ t\left(\EE_{\nu}\left[e^{-\tilde{c}_{nm}}\right]\right)^{1- t}$$

\no The monotone limit theorem enables us to take the limit with
relation to m and then to take it again with respect to n to get the
result.\nqed

\section{\bf{Extension of the map $u\mapsto W^u$}\label{ex}}

\subsection{\bf{Extension of the map $u\mapsto W^u$ for invertibility
    results}}

\no Invertibility results can give stochastic differential equations
solutions in certain cases, so it can be useful to extend these
results to a larger domain.

\begin{definition}
\label{DD}
Set $\DD$ a subset of $G_0(\nu,\beta)$ such that the map $u\mapsto W^u$ can be extended to $\DD$ while verifying the following conditions:\newline
\no (i) $\D\subset\DD\subset G_0(\nu,\beta)$\newline
\no (ii) For every $u\in\DD$, $W^u$ is adapted.\newline
\no (iii)  For every $u\in\DD$, the law of $W^u$
under $\tilde{\nu}^u$ is the same as the law of W under $\nu$, where
$\tilde{\nu}^u$ is defined by
$\frac{d\tilde{\nu}}{d\nu}=\rho(-\delta_\beta\ u)$\newline
\no (iv) For every $u\in\DD$,
$$\beta\circ W^u=\beta+u$$\newline
\no (v) For every $u,v\in\DD$ such that $v+u\circ W^v\in\DD$
$$W^u\circ W^v=W^{v+u\circ W^v}\;\nu-a.s.$$\newline
\no (vi) There exists $\DDD$ such that $\DD\subset \DDD\subset
L^0_a(\nu,H)$, $\DD=\DDD\cap G_0(\nu,\beta)$ and for every $u\in\DD$ such that the equation $u+v\circ W^u$ has
a solution in $G_0(\nu,\beta)$, this equation has a solution in
$\DDD$.
\end{definition}

\begin{proposition}
\no Set $u\in\DD$, for every bounded measurable $f:W\rightarrow\RR$
$$\EE_\nu\left[f\right]=\EE_\nu\left[f\circ W^u\rho(-\delta_\beta u)\right]$$
Moreover
$$W^u\nu\sim\nu$$
\end{proposition}

\nproof The proof is the same as the case $u\in\D$, see section \ref{fr}.\nqed

\no \begin{remark} $\D$ verify the set of conditions above.
\end{remark}\newline

\begin{theorem}
\label{extd2b}
\no For every $u\in \DD\cap L^2(\nu,H)$, we have
$$H(W^u\nu|\nu)\leq\frac{1}{2}\EE_\nu\left[|u|_H^2\right]$$
\no and the three following propositions are equivalent:\newline
\no 1)$H(W^u\nu|\nu)=\frac{1}{2}\EE_\nu\left[|u|_H^2\right]$\newline
\no 2)there exists $v\in\DD$ such that
\beaa W^u\circ W^v=W^v\circ W^u=I_\W\;\nu-a.s.\eeaa
\beaa \frac{dW^u\nu}{d\nu}&=&\rho(-\delta_\beta v)\\
\frac{dW^v\nu}{d\nu}&=&\rho(-\delta_\beta u)\eeaa
\no 3) $W^u$ is $\nu$-a.s. left-invertible
\end{theorem}

\nproof Set $u\in  \D$, we prove that $H(W^u\nu|\nu)\leq\frac{1}{2}\EE_\nu\left[|u|_H^2\right]$ exactly as in proposition \ref{i1}\newline

\no For the second assertion, we prove $(1)\Rightarrow (2)$ first. Exactly as in the proof of theorem \ref{inv}, we obtain

\begin{equation}\label{2beg}L\circ  W^u\rho(-\delta_{\beta} u)= 1\end{equation}

\no Since $\beta$ is a $\nu$-Brownian motion, there exists $v\in
G_0(\nu,\beta)$ such that
$L=\rho(-\delta_{\beta} v)$.
\newline
\no We apply the exponential to \ref{2beg} to get:
\beaa 0&=&\delta_{\beta} v\circ  W^u+\frac{1}{2}|v\circ W^u|^2_H+\delta_{\beta}
u+|u|^2_H\eeaa
\no We have:
\beaa \delta_{\beta} v\circ  W^u &=&\int_0^1\dot{v}(s)\circ W^ud\beta(s) +\langle v\circ  W^u,u\rangle_H\eeaa

\no so finally we have:
\begin{equation}\label{1atest} 0= \delta_{\beta}(v\circ W^u+u)+\frac{1}{2}|v\circ  W^u+u|_H^2\end{equation}
According to Girsanov theorem $\beta+v$ is a $ W^u\nu$-Brownian motion,
so:
\beaa \EE_{\nu}\left[L\log L\right] &=& \EE_{ W^u\nu}\left[\log
    L\right]\\
&=&\EE_{ W^u\nu}\left[-\int_0^1\dot{v}(s)d\beta(s)-\frac{1}{2}\int_0^1|\dot{v}(s)|^2ds\right]\\
&=&\frac{1}{2}\EE_{ W^u\nu}\left[\int_0^1|\dot{v}(s)|^2ds\right]\\
&=&\frac{1}{2}\EE_{\nu}\left[|v\circ  W^u|_H^2\right]\eeaa

\no So $v\circ  W^u\in L^2_a(\nu,H)$ and we can take the expectation
with respect to $\nu$ in \ref{1atest} to obtain $u+v\circ  W^u=0$
$\nu$-a.s. \no Condition (vi) gives the existence of $\tilde{v}\in\DDD$ such that $\nu$-a.s.
$$u+\tilde{v}\circ W^u=0$$
\no We have $v=\tilde{v}$ $W^u\nu$-a.s. so $v=\tilde{v}$ $\nu$-a.s. since $W^u\nu\sim\nu$, which implies $v\in\DD$ and condition (iv) gives $\nu$-a.s.
$$W^v\circ W^u=I_\W$$
\no Proposition \ref{lrinv} gives $\nu$-a.s.
$$W^u\circ W^v=I_\W$$
\no Finally, set $f\in C_b(\W)$,
\beaa \EE_\nu\left[f\circ W^u\right]&=&\EE_\nu\left[f\circ W^u\circ
  W^v\rho(-\delta_\beta v)\right]\\
&=&\EE_\nu\left[f\rho(-\delta_\beta v)\right]\\
\EE_\nu\left[f\circ W^v\right]&=&\EE_\nu\left[f\circ W^v\circ
  W^u\rho(-\delta_\beta u)\right]\\
&=&\EE_\nu\left[f\rho(-\delta_\beta u)\right]\eeaa
\no which gives
\beaa \frac{dW^u\nu}{d\nu}&=&\rho(-\delta_\beta v)\\
\frac{dW^v\nu}{d\nu}&=&\rho(-\delta_\beta u)\eeaa

 \no $(2)\Rightarrow (3)$ is immediate. Now we prove $(3)\Rightarrow (1)$. We still denote $L=\frac{dW^u\nu}{d\nu}$.

\no Assume that $ W^u$ admits a left inverse V. Set
$v=-u\circ V$.\newline
\no We have $\nu$-a.s.
$$v\circ W^u=-u$$
\no and
$$\EE_{ W^u\nu}\left[1_{\int_0^1|\dot{v}(s)|^2ds<\infty}\right]=\EE_{\nu}\left[1_{\int_0^1|\dot{u}(s)|^2ds<\infty}\right]=1$$
\no so $v\in L^0( W^u\nu,H)$ and $v\in L^0(\nu,H)$ since $ W^u\nu\sim\nu$.\newline
\no Now set $\dot{v}^n=\max(n,\min(\dot{v},-n))$, $\dot{v}^n\circ
 W^u$ is adapted. Set $A\in L^2(dt\times d\nu)$ an adapted process, we
have:
\beaa\EE_{\nu}\left[\rho(-\delta_{\beta} u)\int_0^1\dot{v}^n(s)\circ W^u
  A(s)\circ  W^uds\right]&=&\EE_{\nu}\left[\int_0^1\dot{v}^n(s)
  A(s)ds\right]\\
&=&\EE_{\nu}\left[\int_0^1\EE_{\nu}\left[\left.\dot{v}^n(s)\right|\F(s)\right]
  A(s)ds\right]\\
&=&\EE_{\nu}\left[\rho(-\delta_{\beta} u)\int_0^1\EE_{\nu}\left[\left.\dot{v}^n(s)\right|\F(s)\right]\circ W^u
  A(s)\circ  W^uds\right]\eeaa
\no So $\EE_{\nu}\left[\left.\dot{v}^n(s)\right|\F(s)\right]\circ W^u=\dot{v}^n(s)\circ
 W^u$ $ds\times d\nu$-a.s. which implies
$\EE_{\nu}\left[\left.\dot{v}^n(s)\right|\F(s)\right]=\dot{v}^n(s)$
$ds\times d\nu$-a.s. since $ W^u\nu\sim\nu$.

\no An algebraic calculation gives $\nu$-a.s.
$$\rho(-\delta_{\beta} v)\circ  W^u\rho(-\delta_{\beta} u)=1$$
\no Now set $g\in C_b(W,\RR_+)$, we have:
\beaa \EE_{\nu}\left[gL\right]&=&\EE_{\nu}\left[g\circ  W^u\right]\\
&=&\EE_{\nu}\left[g\circ
   W^u\rho(-\delta_{\beta} v)\circ
   W^u\rho(-\delta_{\beta} u)\right]\\
&\leq&\EE_{\nu}\left[g\rho(-\delta_{\beta} v) \right]\eeaa

\no So $L\leq \rho(-\delta_\beta v)$ $\nu$-a.s. and since $\EE_\nu\left[\rho(-\delta_\beta v)\right]=1$, we have
$$L\circ  W^u\rho(-\delta_{\beta} u)=1$$
\no and we can compute $H( W^u\nu|\nu)$:
\beaa H( W^u\nu|\nu)&=&\EE_{\nu}[L\log L]\\
&=&\EE_{\nu}[\log L\circ  W^u]\\
&=&\EE_{\nu}\left[-\log\rho(-\delta_{\beta} u)\right]\\
&=&\frac{1}{2}\EE_{\nu}[|u_H^2|]\eeaa
\nqed

\no As in section \ref{si}, the proof of theorem \ref{extd2b} give the following additional results.

\begin{corollary}
Set $u\in\DD$, we have
$$L\circ W^u\EE_\nu\left[\left.\rho(-\delta_\beta
    u)\right|\F_1^{W^u}\right]=1$$
\no and if $W^u$ is $\nu$-a.s. left-invertible, we have
$$L\circ W^u\rho(-\delta_\beta
    u)=1$$
\end{corollary}

\no In certain cases invertibility results lead to the existence of a strong solutions of stochastic differential equations.

\begin{theorem}
\label{sde}
 Assume that for every $u\in\DD$, we can write $\nu$-a.s.
$$W^u=I_\W+w^u$$
\no with $w^u\in L^0_a(\nu,H)$.\newline
\no Set $u\in\DD$,
$$H(W^u\nu|\nu)=\frac{1}{2}\EE_\nu\left[|u|_H^2\right]$$
\no if and only if there exists $v\in\DD$ such that $W^v$ is a strong solution to
$$dW^v(t)=-\dot{\overbrace{w^u}}(t)\circ W^vdt+dW(t)$$

\end{theorem}

\nproof Assume that $H(W^u\nu|\nu)=\frac{1}{2}\EE_\nu\left[|u|_H^2\right]$ , according to theorem \ref{extd2b}, there exists $v\in\DD$ such that $\nu$-a.s.
$$W^v\circ W^u=W^u\circ W^v=W$$
\no Since $W^u=W+w^u$, we have
$$W^v+w^u\circ W^v=W$$
\no and $W^v$ is a strong solution to
$$dW^v(t)=-\dot{\overbrace{w^u}}(t)\circ W^vdt+dW(t)$$\newline
\no Conversely, assume the existence of $v\in\DD$ such $W^v$ is a strong solution to
$$dW^v(t)=-\dot{\overbrace{w^u}}(t)\circ W^vdt+dW(t)$$
\no We have $W^u\circ W^v=I_W$ $\nu$-a.s., and $W^v\circ W^u=I_W$ $W^v\nu$-a.s. Since $W^v\nu\sim\nu$, we can conclude with theorem \ref{extd2b}.\nqed

\subsection{\bf{Extension of the map $u\mapsto W^u$ for variational calculus}}

\begin{theorem}
\label{extd}
For every measurable function $f:\W\rightarrow\RR$ such that $\EE_\nu\left[(|f|+1)e^{-f}\right]<\infty$ and
$$-\log\EE_\nu\left[e^{-f}\right]=\inf_{u\in\D^i}\EE_\nu\left[f\circ
  W^u+\frac{1}{2}|u|_H^2\right]$$
\no we have
$$-\log\EE_\nu\left[e^{-f}\right]=\inf_{u\in\DD\cap L^2_a(\nu,H)}\EE_\nu\left[f\circ W^u+\frac{1}{2}|u|_H^2\right]$$
\end{theorem}

\nproof For $u\in\DD\cap L^2_a(\nu,H)$, we have
\beaa
-\log\EE_\nu\left[e^{-f}\right]&\leq&\EE_{W^u\nu}\left[f\right]+H(W^u\nu|\nu)\\
&\leq&\EE_\nu\left[f\circ W^u+\frac{1}{2}|u|_H^2\right]\eeaa
\no So
\beaa
-\log\EE_\nu\left[e^{-f}\right]&\leq&\inf_{u\in\DD\cap L^2_a(\nu,H)}\EE_\nu\left[f\circ
  W^u+\frac{1}{2}|u|_H^2\right]\\
&\leq&\inf_{u\in\D^i}\EE_\nu\left[f\circ
  W^u+\frac{1}{2}|u|_H^2\right]\eeaa\nqed

\begin{theorem}
 Set $f:\W\rightarrow\RR$ a measurable function
verifying $\EE_{\nu}\left[|f|(1+e^{-f})\right]<\infty$, then if there
exists some $u\in \DD\cap L^2_a(\nu,H)$ such that $ W^u$ is
$\nu$-a.s. left-invertible and
$\frac{d W^u\nu}{d\nu}=\frac{e^{-f}}{\EE_{\nu}\left[e^{-f}\right]}$, then
we have
$$-\log\EE_{\nu}\left[e^{-f}\right]=\inf_{u\in \DD\cap L^2_a(\nu,H)}\EE_{\nu}\left[f\circ  W^u+\frac{1}{2}|u|_H^2\right]$$
\end{theorem}

\nproof Since $ W^u$ is $\nu$-a.s. left invertible and that
$\frac{d W^u\nu}{d\nu}=\frac{e^{-f}}{\EE_{\nu}\left[e^{-f}\right]}$. We have
$$\frac{1}{2}\EE_{\nu}\left[|u|_H^2\right]=H( W^u\nu|\nu)=\EE_{\nu}\left[\frac{e^{-f}}{\EE_{\nu}\left[e^{-f}\right]}\log
  \left(\frac{e^{-f}}{\EE_{\nu}\left[e^{-f}\right]}\right)\right]$$
\no and
\beaa\EE_{\nu}\left[f\circ
   W^u+\frac{1}{2}|u|_H^2\right]&=&\EE_{\nu}\left[\frac{e^{-f}}{\EE_{\nu}\left[e^{-f}\right]}f+\frac{e^{-f}}{\EE_{\nu}\left[e^{-f}\right]}\log\left(\frac{e^{-f}}{\EE_{\nu}\left[e^{-f}\right]}\right)\right]\\
&=&-\log\EE_{\nu}\left[e^{-f}\right]\eeaa
\no and we conclude the proof with last proposition.\nqed

\begin{theorem}
Set $f:\W\rightarrow \RR$  a measurable function such that
$$-\log\EE_{\nu}\left[e^{-f}\right] = \inf_{u\in
  \DD\cap L^2_a(\nu,H)}\EE_\nu\left[f\circ  W^u+\frac{1}{2}|u|_H^2\right]$$
\no Denote this infimum $J_*$. It is attained at $u\in \DD\cap L^2_a(\nu,H)$ if and only
if $ W^u$ is $\nu$-a.s. left-invertible and
$\frac{d W^u\nu}{d\nu}=\frac{e^{-f}}{\EE_{\nu}\left[e^{-f}\right]}$.\newline
\end{theorem}

\nproof The direct implication is given by last theorem.
\no Conversely, if $ W^u$ is not $\nu$-a.s. left-invertible,
$H( W^u\nu|\nu)<\frac{1}{2}\EE_{\nu}\left[|u|_H^2\right]$ and
\beaa -\log\EE_{\nu}\left[e^{-f}\right]=\inf_{\theta\in
  \mathcal{P}}\left(\EE_\theta[f]+H(\theta|\nu)\right)&\leq& \inf_{\alpha\in
  \DD\cap L^2_a(\nu,H)}\EE_{ W^\alpha \nu}\left[f\right]+H( W^\alpha \nu|\nu)\\
&\leq&\EE_{ W^u \nu}\left[f\right]+H( W^u \nu|\nu)\\
&<&\EE_{\nu}\left[f\circ  W^u+\frac{1}{2}|u|_H^2\right]\eeaa
\no which is a contradiction.\newline
\no We get $\frac{d W^u\nu}{d\nu}=L$ by uniqueness of the minimizing
measure of $\inf_{\theta\in
  \mathcal{P}}\left(\EE_\theta[f]+H(\theta|\nu)\right)$.\newline
\nqed

\section{\bf{Examples}}

\no In this section we discuss several examples that fit into the framework we elaborated. Each time, we prove that the conditions of section \ref{fr} and definition \ref{DD} are satisfied. This ensures that every result from section \ref{fr} to \ref{ex} apply, except theorems \ref{pl} and \ref{sde} which require additional hypothesis. We also discuss whether these theorems apply or not.

\subsection{\bf{Diffusion}}

\no Set $m\leq d\in \NN^*$ such that $m+d=n$, $c\in\RR$, $\sigma:\RR^m\rightarrow\M_{m,d}(\RR)$ and $b:\RR^m\rightarrow\RR^m$ bounded and lipschitz functions. $\sigma_i$ will denote the
i-th column of $\sigma$. Notice that every matrix
will be identified with its canonical linear operator. Set $(\Omega,\theta,(\mathcal{G}_t))$ a probability space, V a $\theta$-Brownian motion on
$\Omega$ with values in $\RR^d$. Set Y a $\RR^m$-valued strong solution of the
stochastic differential equation:
$$Y(t)=c+\int_0^t\sigma(Y(s))dV(s)+\int_0^tb(Y(s))ds$$
\no on $(\Omega,\theta,(\mathcal{G}_t),B)$. The hypothesis on $\sigma$ and b ensure the
existence and uniqueness of Y if we impose its paths to be continuous.\newline
\no We denote $\mu$
the Wiener measure on $C([0,1],\RR^d)$ and $\mu^X$ the image measure of Y.\newline\no We define the processes X and B on $\W$ by:
\beaa X(t)&:&(w,w')\in \W\mapsto w(t)\in\RR^m\\
B(t)&:&(w,w')\in \W\mapsto w'(t)\in\RR^d\eeaa

\begin{proposition}
Under $\mu^X\times\mu$, the law of X is $\mu^X$, B is a Brownian
motion and they are independent. There exists $\theta,\eta$ such that if we define $\beta_\grandX$ as
$$\beta_\grandX=\int_0^.\theta(X(s))dM(s)+\int_0^.\eta(X(s))dB(s)$$
\no $\beta_\grandX$ is a $\mu^X\times\mu$-Brownian motion and
$\mu^X\times\mu$-a.s.
$$X=c+\int_0^.\sigma(X(s))d\beta_\grandX(s)+\int_0^.b(X(s))ds$$
\end{proposition}
\nproof See \cite{a1}.\nqed

\no This construction of $\beta_\grandX$ is taken from \cite{ro}.

\begin{definition} We denote
$$\grandX=(X,\beta_\grandX)$$
\no and $\loigrandX$ its image measure.\newline
X is a $\loigrandX$ path-continuous strong solution of the
stochastic differential equation
\beaa X=&c+\int_0^.\sigma(X(s))d\beta_\grandX(s)+\int_0^.b(X(s))ds\eeaa
\no For $u\in G_0(\loigrandX,\beta_X)$, set $\beta_\grandX^u=\beta+u$ and $X^u$ the $\loigrandX$-a.s. path-continuous strong solution
of the stochastic differential equation
$$X^u=c+\int_0^.\sigma(X^u(s))d\beta_\grandX^u(s)+\int_0^.b(X^u(s))ds$$
\no Finally, we denote
\beaa\grandX^u&=&(X^u,\beta_\grandX+u)\eeaa
\end{definition}

\begin{theorem} $\left(\W,\loigrandX,\beta_\grandX,(\grandX^u)_{u\in\D}\right)$ verify the
    conditions of section \ref{fr}. $\left(\W,\mu_a,\beta_\grandX,\left(\grandX^u\right)_{u\in G_0(\loigrandX,\beta_\grandX)}\right)$
    verify the conditions of definition \ref{DD}.
\end{theorem}

\nproof See \cite{a1}.\nqed

\begin{corollary} It is clear that for every $u\in\D$, we clearly have $\loigrandX$-a.s.
$$\grandX^u(w)=\grandX^{u(w)}(w)$$
\no so theorem \ref{pl} applies.
\end{corollary}

\begin{corollary} Assume that $\sigma\in\RR$, then theorem \ref{sde} applies. Set $u\in G_2(\loigrandX,\beta_\grandX)$ and denote
$$\dot{\overbrace{w_\grandX^u}}(t)=\left(\dot{u}(t)+b(X^u(t))-b(X(t)),\dot{u}(t)\right)$$
\no We have $\loigrandX$-a.s.
$$\grandX^u=I_\W+w_\grandX^u$$
\no so
$$H(\grandX^u\loigrandX|\loigrandX)=\frac{1}{2}\EE_{\loigrandX}\left[|u|_H^2\right]$$
\no if and only if there exists $v\in G_0(\loigrandX,\beta_\grandX)$ such that $\grandX^v$ is a strong solution to the stochastic differential equation:
$$d\grandX^v(t)=-\dot{\overbrace{w_\grandX^u}}(t)dt\circ \grandX^v+dW(t)$$
\end{corollary}

\subsection{\bf{Brownian bridge}}
\no We still denote $\mu$ the Wiener measure on $\W$. Set $a\in\R^n$, we denote $\mu_a$ the measure on $\W$ such that for any bounded measurable function f we have
$$\EE_{\mu_a}\left[f\right]=\EE_\mu\left[f|W_1=a\right]$$
\no $\mu_a$ can also be defined as follow : let $\mathcal{E}_a$ be the Dirac measure in a, $\mathcal{E}_a(W_1)$ is a positive Wiener distributions hence it defines a Radon measure $\nu_a$ on $\W$, then
$$\mu_a=\left(\frac{1}{2\pi}\right)^n\nu_a$$

\no We recall the definition of a Brownian bridge:
\begin{definition}
Set $(\Omega,\mathcal{G},Q)$ a probability space. An a-Brownian bridge X under a probability Q is a path-continuous Gaussian process such
that $\EE_Q\left[X(t)\right]=at$ and $cov(X(s),X(t))=\left((s\wedge
  t)-st\right)I_d$
\end{definition}

\begin{proposition}
W is an a-Brownian bridge under $\mu_a$, and the process $\beta_a$ defined as
$$\beta_a(t)=W(t)-at+\int_0^t\frac{W(s)-as}{1-s}ds$$
\no is a Brownian motion under $\mu_a$ and the filtrations of $\beta_a$ and
W completed with respect to $\mu_a$ are equal. Moreover, we have
$$W(t)=at+(1-t)\int_0^t\frac{d\beta_a(s)}{1-s}$$
\end{proposition}

\nproof It is easy to verify that the process Z given by
$Z(t)=W(t)+at-tW_1$ is an a-Brownian bridge under $\mu$ and is independent of
$W_1$. If f is a bounded continuous function on W, we have
\beaa \EE_\mu\left[f(Z)\right]&=&\EE_\mu\left[f(Z)|W_1=a\right]\\
&=&\EE_\mu\left[f(W)|W_1=a\right]\\
&=&\EE_\mu\left[f|W_1=a\right]\\
&=&\EE_{\mu_a}\left[f\right]\\
&=&\EE_{\mu_a}\left[f(W)\right]\eeaa
\no So W is indeed an a-Brownian bridge under $\mu_a$.\newline
\no Now consider the process X given by
$$X(t)=(1-t)W\left(\frac{t}{1-t}\right)+at$$
\no It is easy again to verify that X is a Brownian bridge under
$\mu$. Now consider
$$M=\int_0^.\frac{dW(s)}{1-s}$$
\no M is a continuous martingale under $\mu$ and
$$\left\langle
  M_i,M_j\right\rangle(t)=\delta_{ij}\int_0^t\frac{ds}{(1-s)^2}=\delta_{ij}\frac{t}{1-t}$$
\no So the Dubins-Schwartz theorem ensures that M and
$W(\frac{.}{1-.})$ have the same distribution under $\mu$, so X and
$\tilde{X}$ have the same distribution, where $\tilde{X}$ is given by
\beaa \tilde{X}(t)&=&(1-t)\int_0^t\frac{dW(s)}{1-s}+at\\
&=&W(t)-\int_0^t\frac{\tilde{X}(s)-as}{1-s}ds+at\eeaa
\no the last equality coming from Ito's formula.\newline
\no W is a Brownian motion under $\mu$ and the law of
$(\tilde{X},W)$ under $\mu$ is the same as the law of
$(W,\beta_a)$ under $\mu_a$ so $\left(\beta_a(t),t\in [0,1)\right)$ is a Brownian motion under
$\mu_a$.\newline
\no We recall that the filtrations we are considering are all
completed with respect to $\mu_a$.  From the expression of $\beta_a$, we have clearly for any $t\in
[0,1)$,
$$\F_t^{\beta_a}\subset\F_t^W$$
\no Furthermore, $s\mapsto\frac{1}{1-s}$ being lipschitz on any
$[0,t]$ with $t<1$, W is the strong solution of a stochastic differential equation relative to $\beta_a$ and we have for any $t<1$,
$$\F_t^W\subset \F_t^{\beta_a}$$
\no W being $\mu_a$-a.s. path continuous, we have
$$\bigcup_{t<1}\F_t^W=\F_1^W$$
\no On the other hand, since $\left(\beta_a(t),t\in [0,1)\right)$ is a Brownian motion, $\left(\left\langle\beta_a\right\rangle(t),t\in [0,1)\right)$ is bounded by 1, so $\beta_a(t)$ converges $\mu_a$-a.s. and in $L^2(\mu_a,\RR^n)$. Denote $\beta_a(1)$ the limit, for $s\in[0,1)$ we have
$$Cov(\beta_a(s),\beta_a(t))=\lim_{t'\rightarrow 1}Cov(\beta_a(s),\beta_a(t'))=s\wedge t'$$
\no So
$$\bigcup_{t<1}\F_t^{\beta_a}=\F_1^{\beta_a}$$
\no $(\beta_a(t),t\in [0,1])$ is a $\mu_a$-Brownian motion and $\beta_a$ and W have the same filtration.\nqed

\no The following remark will be useful in next section.\newline
\begin{remark}
For $a\in\RR^n$ and $t\in [0,1]$, we have $\mu_a$-a.s.
$$\beta_a(t)=W_t+\int_0^t\frac{W_s-a}{1-s}ds$$
\end{remark}

\begin{definition}
\no For $u\in G_0(\mu_a,\beta_a)$, we denote $\beta_a^u=\beta_a+u$.
\end{definition}

\begin{proposition}
\no Set $u\in G_0(\mu_a,\beta_a)$, then there exists a unique $\mu_a$-a.s. path
continuous process $W_a^u$ such that
$$W_a^u(t)=\beta_a^u(t)+at-\int_0^t\frac{W_a^u(s)-as}{1-s}ds$$
\no Furthermore, we have
\beaa W_a^u(t)&=&at+ (1-t)\int_0^t\frac{d\beta_a^u(s)}{1-s}\\
&=&W(t)+\int_0^t\left(\dot{u}(s)-\int_0^s\frac{\dot{u}(r)}{1-r}dr\right)ds\eeaa
\end{proposition}

\nproof Set $u\in G_0(\mu_a,\beta_a)$, for $t<1$, straight calculation gives
$$at+(1-t)\int_0^t\frac{d\beta_a^u(s)}{1-s}=W+\int_0^t\left(\dot{u}(s)-\int_0^s\frac{\dot{u}(r)}{1-r}dr\right)ds$$
\no Define $W_a^u$ on $[0,1)$ as
$$W_a^u(t)=at+(1-t)\int_0^t\frac{d\beta_a^u(s)}{1-s}$$
\no The Ito formula gives
$$W_a^u(t)=\beta_a^u(t)+at-\int_0^t\frac{W_a^u(s)-as}{1-s}ds$$
\no $x\mapsto\frac{1}{1-x}$ being lipschitz on every $[0,t]$ for $t<1$,
the $\mu_a$-a.s. pathwise uniqueness is true on every $[0,t]$, hence on
$[0,1)$.\newline
\no It remains to prove that there is no explosion in 1. Set
$\tilde{\mu}_a^u$ the measure on $\W$ defined by
$$\frac{d\tilde{\mu}_a^u}{d\mu_a}=\rho(-\delta_{\beta_a} u)$$
\no Since $u\in G_0(\mu_a,\beta_a)$, it is clear that the law of $W_a^u$ under $\tilde{\mu}_a^u$ is the same as the law
of $W_a$ under $\mu_a$ so
$$\tilde{\mu}_a^u\left(\limsup_{t\rightarrow
    1}\left|W_a^u(t)\right|=\infty\right)=\mu_a\left(\limsup_{t\rightarrow
    1}\left|W_a(t)\right|=\infty\right)=0$$
\no $\tilde{\mu}_a^u\sim\mu_a$ so
$$\mu_a\left(\limsup_{t\rightarrow
    1}\left|W_a^u(t)\right|=\infty\right)=0$$\nqed

\begin{theorem} $\left(\W,\mu_a,\beta_a,(W_a^u)_{u\in\D}\right)$ verify the
    conditions of section \ref{fr}. $\left(\W,\mu_a,\beta_a,\left(W_a^u\right)_{u\in G_0(\mu_a,\beta_a)}\right)$
    verify the conditions of definition \ref{DD}.
\end{theorem}
\nproof We have $W_a^0=W$ and $\beta_a$ is a $\mu_a$-Brownian motion. Now we
just have to verify the conditions of definition \ref{DD} since those
imply conditions (iii) to (v) of section \ref{fr}.\newline
\no (i), (ii)  and (iii) are clear, so is (vi) taking $\DDD=G_0$.\newline
\no Set $u\in G_0(\mu_a,\beta_a)$, we have
 \beaa \beta_a\circ
W_a^u(t)&=&\left(W(t)-at+\int_0^t\frac{W(s)-as}{1-s}ds\right)\circ W_a^u\\
&=&W_a^u(t)-at+\int_0^t\frac{W_a^u(s)-as}{1-s}ds\\
&=&\beta_a^u(t)\eeaa
\no so condition (iv) is verified. Now set $v\in G_0(\mu_a,\beta_a)$ such that $v+u\circ W^v\in G_0(\mu_a,\beta_a)$, we have
\beaa W_a^u(t)\circ
W_a^v&=&\left(W(t)+\int_0^t\left(\dot{u}(s)-\int_0^s\frac{\dot{u}(r)}{1-r}dr\right)ds\right)\circ W_a^v\\
&=&W_a^v(t)+\int_0^t\left(\dot{u}(s)\circ W_a^v-\int_0^s\frac{\dot{u}(r)\circ
    W_a^v}{1-r}dr\right)ds\\
&=&W(t)+\int_0^t\left(\dot{v}(s)+\dot{u}(s)\circ W^v_a-\int_0^s\frac{\dot{v}(r)+\dot{u}(r)\circ
 W^v(r)}{1-r}dr\right)ds\\
&=&W^{v+u\circ W_a^v}\eeaa
\no which gives condition (v).\nqed

\begin{corollary} It is clear that for every $u\in\D$, we clearly have $\mu_a$-a.s.
$$W_a^u(w)=W_a^{u(w)}(w)$$
\no so theorem \ref{pl} applies.
\end{corollary}

\begin{corollary}
Theorem \ref{sde} applies. Set $u\in G_2(\mu_a,\beta_a)$, we have $\mu_a$-a.s.
$$W_a^u=I_\W+\int_0^.\dot{u}(t)-\int_0^t\frac{\dot{u}(s)}{1-s}dsdt$$
\no so
$$H(W_a^u\mu_a|\mu_a)=\frac{1}{2}\EE_{\mu_a}\left[|u|_H^2\right]$$
\no if and only if there exists $v\in G_0(\mu_a,\beta_a)$ such that $W_a^v$ is a strong solution to the stochastic differential equation:
$$dW_a^v(t)=-\left(\dot{u}(t)-\int_0^t\frac{\dot{u}(s)}{1-s}ds\right)dt\circ W_a^v+dW(t)$$
\end{corollary}

\subsection{\bf{Loop measure}}

\no We keep the notations of last section. Denote
$$S=\left\{a\in\RR^n,|a|=1\right\}$$
\no and set $\alpha:S \rightarrow \RR_+$ a locally lipschitz
function such that $\left\{x,\alpha(x)\neq 0\right\}$ is of strictly
positive measure for the Lebesgue measure on S and
$$\int_S\alpha(a)da=1$$
\no We define the measure $\nu_l$ as follow:
for any bounded measurable function f on W, we
set
$$\EE_{\nu_l}[f]=\int_S\alpha(a)\EE_{\mu_a}[f]da$$
\no For more on loop measures, see Fang's work in \cite{fan}.
\begin{definition}
We denote
\beaa h_a:(t,x)&\in&
[0,1)\times\RR^n\mapsto\left(\frac{1}{\pi(1-t)}\right)^{\frac{n}{2}}\exp\left(\frac{-\left|x-a\right|^2}{2(1-t)}\right)\\
h:(t,x)&\in& [0,1)\times\RR^n\mapsto\int_S \alpha(a)h_a(t,x)da\eeaa
\end{definition}

\begin{proposition}
\label{dl}
Set $a\in\RR^n$ and $t\in[0,1)$, then
$$\left.\frac{d\mu_a}{d\mu}\right|_{\F^W_t}=h_a(t,W_t)$$
\end{proposition}

\nproof For convenience we consider the case $n=1$, the general proof
is the same. Every $\F^W_t$ measurable $f:\W\rightarrow\RR$ is a function of
$W(.\wedge t)$, hence of $\beta_a(.\wedge t)$. $\beta_a(.\wedge t)$ is a Brownian motion on $[0,t]$
under $\mu_a$. Now denote $\tilde{\mu}$ the probability measure on $\W$ given by
$$\frac{d\tilde{\mu}}{d\mu}=\exp\left(-\int_0^t\frac{W(s)-a}{1-s}dW(s)-\frac{1}{2}\int_0^t\left(\frac{W(s)-a}{1-s}\right)^2ds\right)$$
\no According to Girsanov theorem, $\beta_a(.\wedge t)$ is also a Brownian motion under $\tilde{\mu}$ and
\beaa
\left.\frac{d\mu_a}{d\mu}\right|_{\F^W_t}&=&\left.\frac{d\tilde{\mu}}{d\mu}\right|_{\F^W_t}\\
&=&\exp\left(-\int_0^t\frac{W(s)-a}{1-s}dW(s)-\frac{1}{2}\int_0^t\left(\frac{W(s)-a}{1-s}\right)^2ds\right)\eeaa
\no Finally, Ito formula gives
$$\left(\frac{1}{\pi(1-t)}\right)^{\frac{1}{2}}\exp\left(\frac{-\left|W(s)-a\right|^2}{2(1-t)}\right)=\exp\left(-\int_0^t\frac{W(s)-a}{1-s}dW(s)-\frac{1}{2}\int_0^t\left(\frac{W(s)-a}{1-s}\right)^2ds\right)$$\nqed

\begin{proposition}
Set $t\in [0,1)$, we have
$$\left.\frac{d\nu}{d\mu}\right|_{\F_t^W}=h(t,W_t)$$
\end{proposition}

\nproof Set $C\in\F_t^W$, Fubini-Tonelli theorem gives
\beaa \EE_\nu\left[1_C\right]&=& \int_S\alpha(a)\EE_{\mu_a}\left[1_C\right]da\\
&=&\int_S\alpha(a)\EE_\mu\left[1_C h_a(t,W(t))\right]da\\
&=&\EE_\mu\left[1_C\int_S \alpha(a)h_a(t,W(t))da\right]\\
&=&\EE_\mu\left[1_Ch(t,W(t))\right]\eeaa\nqed

\begin{proposition}
Define
$$\beta_{\lo}(t)= W(t)-\int_0^t\frac{h'(s,W(s))}{h(s,W(s))}ds$$
\no where h' designates the partial derivative of h with respect to
x.\newline
\no Then $\beta_{lo}$ is a $\nu_l$ Brownian motion
and the filtration of W and $\beta_{\lo}$ completed with respect to
$\nu_l$ are equal.
\end{proposition}

\nproof Notice that every filtration we consider here is completed with respect to $\nu_l$. The fact that $\left(\beta_{\lo}(t),t\in [0,1)\right)$ is a
$\nu$-Brownian motion is direct consequence of proposition \ref{dl} and the expression of $\beta_{\lo}$ gives that for $t>1$,
$$\F_t^{\beta_{\lo}}\subset\F_t^W$$
\no On the other hand, for $t<1$ since $s\mapsto\frac{h'(s,x)}{h(s,x)}$ is lipschitz on $[0,t]$ and $x\mapsto\frac{h'(s,x)}{h(s,x)}$ is lipschitz on $\{x\in\R^n,|x|\leq k\}$ for any $k>0$ so W is the strong solution of a stochastic differential equation relative to $\beta_{lo}$ and
$$\F_t^W\subset\F_t^{\beta_{\lo}}$$
\no W being $\mu_a$-a.s. path continuous, we have
$$\bigcup_{t<1}\F_t^W=\F_1^W$$
\no On the other hand, since $\left(\beta_{lo}(t),t\in [0,1)\right)$ is a Brownian motion, $\left(\left\langle\beta_{lo}\right\rangle(t),t\in [0,1)\right)$ is bounded by 1, so $\beta_{lo}(t)$ converges $\mu_{lo}$-a.s. and in $L^2(\nu_l,\RR^n)$. Denote $\beta_{lo}(1)$ the limit, for $s\in[0,1)$ we have
$$Cov(\beta_{lo}(s),\beta_{lo}(t))=\lim_{t'\rightarrow 1}Cov(\beta_{lo}(s),\beta_{lo}(t'))=s\wedge t'$$
\no and finally
$$\bigcup_{t<1}\F_t^{\beta_{lo}}=\F_1^{\beta_{lo}}$$
\no $(\beta_{lo}(t),t\in [0,1])$ is a $\nu_l$-Brownian motion and $\beta_{lo}$ and W have the same filtrations.\nqed

\begin{definition}
\no For $u\in G_0(\nu_l,\beta_{lo})$, we denote $\beta_{\lo}^u=\beta_{\lo}+u$.
\end{definition}

\begin{proposition}
\label{ul}
\no Set $u\in G_0(\nu_l,\beta_{lo})$, then there exists a unique $\nu_l$-a.s. path
continuous process $W_{\lo}^u$ such that
$$W_{\lo}^u(t)=\beta_{\lo}^u(t)+\int_0^t\frac{h'(s,W_{lo}^u(s))}{h(s,W_{lo}^u(s))}ds$$
\end{proposition}

\nproof Set $t<1$, since $s\mapsto\frac{h'(s,x)}{h(s,x)}$ is lipschitz on $[0,t]$ for any $t<1$ and $x\mapsto\frac{h'(s,x)}{h(s,x)}$ is lipschitz on $\{x\in\R^n,|x|\leq k\}$ for any $k>0$, there exists a unique $\nu_l$-a.s. path-continuous process $\left(W_{lo}^u,u\in[0,1)\right)$ such that for any $t<1$,
$$W_{\lo}^u(t)=\beta_{\lo}^u(t)+\int_0^t\frac{h'(s,W_{lo}^u(s))}{h(s,W_{lo}^u(s))}ds$$
\no It remains to prove that there is no explosion in 1. Set
$\tilde{\nu}_l^u$ the measure on $\W$ defined by
$$\frac{d\nu_{l}^u}{d\nu_l}=\rho(-\delta_{\beta_{lo}} u)$$
\no Since $u\in G_0(\nu_l,\beta_{lo})$, it is clear that the law of $W_{\lo}^u$ under $\nu_{l}^u$ is the same as the law
of $W_{\lo}$ under $\nu_l$ so
$$\nu_{l}^u\left(\limsup_{t\rightarrow
    1}\left|W_{\lo}^u(t)\right|=\infty\right)=\nu_l\left(\limsup_{t\rightarrow
    1}\left|W_{\lo}(t)\right|=\infty\right)=0$$
\no $\tilde{\nu}_{l}^u\sim\nu_l$ so
$$\nu_l\left(\limsup_{t\rightarrow
    1}\left|W_{\lo}^u(t)\right|=\infty\right)=0$$\nqed

\begin{theorem} $\left(\W,\nu_l,\beta_{\lo},(W_{\lo}^u)_{u\in\D}\right)$ verify the
    conditions of section \ref{fr}.\newline \no$\left(\W,\nu_l,\beta_{lo},\left(W_{\lo}^u\right)_{u\in G_0(\nu_l,\beta_{lo})}\right)$
    verify the conditions of definition \ref{DD}.
\end{theorem}
\nproof We have $W_{\lo}^0=W$ and $\beta_{\lo}$ is a $\nu_l$-Brownian motion. Now we
just have to verify the conditions of definition \ref{DD} since those
imply conditions (iii) to (v) of section \ref{fr}.\newline
\no (i), (ii) and (iii) are clear, so is (vi) taking $\DDD=G_0$.\newline
\no Set $u\in G_0(\nu_l,\beta_{lo})$, we have
 \beaa \beta_{\lo}\circ
W_{\lo}^u(t)&=&\left(W(t)-\int_0^t\frac{h'(s,W(s))}{h(s,W(s))}ds\right)\circ W_{\lo}^u\\
&=&W_{lo}^u(t)-\int_0^t\frac{h'(s,W_{lo}^u(s))}{h(s,W_{lo}^u(s))}ds\\
&=&\beta_{\lo}^u(t)\eeaa
\no so condition (iv) is verified. Now set $v\in G_0(\nu_l,\beta_{lo})$ such that $v+u\circ W^v\in G_0(\nu_l,\beta_{lo})$, we have
\beaa W_{\lo}^u(t)\circ
W_{\lo}^v&=&\left(\beta_{lo}^u(t)+\int_0^t\frac{h'(s,W(s))}{h(s,W(s))}ds\right)\circ W_{\lo}^v\\
&=&\beta_{lo}(t)+v(t)+u(t)\circ W_{lo}^v +\int_0^t\frac{h'(s,W_{lo}^u(s))}{h(s,W_{lo}^u(s))}ds\eeaa
\no $W^u\circ W^v$ and $W^{v+u\circ W^v}$ are both $\nu_l$-a.s. path continuous so the uniqueness result in proposition \ref{ul} gives $\nu_l$-a.s.
$$W^u\circ W^v=W^{v+u\circ W^v}$$
\no and condition (v) is verified.\nqed

\begin{corollary} It is clear that for every $u\in\D$, we clearly have $\nu_l$-a.s.
$$W_{lo}^u(w)=W_{lo}^{u(w)}(w)$$
\no so theorem \ref{pl} applies.
\end{corollary}

\begin{corollary}
Theorem \ref{sde} applies. Set $u\in G_2(\nu_l,\beta_{lo})$, we have $\nu_l$-a.s.
$$W_{lo}^u=W+\int_0^.\dot{u}(t)+\frac{h'(t,W_{lo}^u(t))}{h(t,W_{lo}^u(t))}-\frac{h'(t,W(t))}{h(t,W(t))}dt$$
\no so
$$H(W_{lo}^u\nu_l|\nu_l)=\frac{1}{2}\EE_{\nu_l}\left[|u|_H^2\right]$$
\no if and only if there exists $v\in G_0(\nu_l,\beta_{lo})$ such that $W_{lo}^v$ is a strong solution to the stochastic differential equation:
$$dW_{lo}^v(t)=-\left(\dot{u}(t)+\frac{h'(t,W_{lo}^u(t))}{h(t,W_{lo}^u(t))}-\frac{h'(t,W(t))}{h(t,W(t))}\right)dt\circ W_{lo}^v+dW(t)$$
\end{corollary}

\subsection{\bf{Diffusing particles without collision}}

\no Set $\sigma,b,\delta,\gamma\in\RR$ such that
$$\sigma^2\leq 2\gamma$$
\no The proof of the following theorem can be found in \cite{shi} or \cite{cl}.
\begin{theorem}
\label{par}
Set $(\Omega,\theta,(\mathcal{G}_t))$ a filtered probability space, $(z_1(0),...,z_n(0))\in\RR^n$ and $B=(B_1,...,B_n)$ a $\RR^n$-valued $\theta$-Brownian motion. We consider the following stochastic differential system:
\beaa Z_1(t)&=&z_1(0)+\sigma B_1(t)+b\int_0^t Z_1(s)ds+ct+\gamma\sum_{j\in\{1,...,n\}\backslash\{1\}}\int_0^t\frac{ds}{Z_1(s)-Z_j(s)}\\
\vdots&&\\
Z_n(t)&=&z_n(0)+\sigma B_n(t)+b \int_0^tZ_n(s)ds+ct+\gamma\sum_{j\in\{1,...,n\}\backslash\{n\}}\int_0^t\frac{ds}{Z_n(s)-Z_j(s)}\eeaa
\no under the condition that $\theta$-a.s. for every $t\in [0,\infty)$
$$Z_1(t)\leq...\leq Z_n(1)$$
\no This system admits a unique strong solution on $(\Omega,\theta,(\mathcal{G}_t),B)$ and the first collision time is $\theta$-a.s. equal to $\infty$.
\end{theorem}

\no Consider $(\Omega,\theta,(\mathcal{G}_t))$ a filtered probability space, $(z_1(0),...,z_n(0))\in\RR^n$ and $B=(B_1,...,B_n)$ a $\RR^n$-valued $\theta$-Brownian motion, and Z the strong solution of the stochastic differential system of theorem \ref{par}. Denote $\nu_{pa}=Z$ the image measure of Z. For $1\leq i\leq n$, denote $W_1,...,W_n$ the coordinates of W and define

$$M_i(t)=W_i(t)-z_i(0)-b\int_0^tW_i(s)ds-ct-\gamma\sum_{j\in\{1,...,n\}\backslash\{i\}}\int_0^t\frac{ds}{W_i(s)-W_j(s)}$$
\no and
$$M=(M_1,...,M_n)$$
\no M is a local martingale and
$$\langle M_i,M_j\rangle(t)=\sigma^2t$$
\no Define
$$\beta_{pa}=\frac{1}{\sigma}M$$
\no Levy theorem clearly ensures that $\beta$ is a $\nu_{pa}$-Brownian motion and we clearly have for every $1\leq i\leq n$,
$$W_i(t)=z_i(0)+\sigma\beta_{pa,i}(t)+b\int_0^tW_i(s)ds+ct+\gamma\sum_{j\in\{1,...,n\}\backslash\{i\}}\int_0^t\frac{ds}{W_i(s)-W_j(s)}$$

\no For $u\in G_0(\nu_{pa},\beta_{pa})$ denote
$$\beta_{pa}^u=\beta_{pa}+u$$
\no and $\nu_{pa}^u$ the probability measure given by
$$\frac{d\nu_{pa}^u}{d\nu_{pa}}=\rho(-\delta_{\beta_pa}u)$$
\no According to Girsanov theorem, $\beta_{pa}+u$ is a Brownian motion under $\nu_{pa}^u$, so according to theorem \ref{par}, there exists a unique $\nu_{pa}^u$-a.s. continuous process $W_{pa}^u=(W_{pa,1}^u,...,W_{pa,n}^u)$ such that $\nu_{pa}^u$-a.s. for every $1\leq i\leq n$
$$W_{pa,i}^u(t)=z_i(0)+\sigma\beta_{pa,i}^u(t)+b\int_0^tW_{pa,i}^u(s)ds+ct+\gamma\sum_{j\in\{1,...,n\}\backslash\{i\}}\int_0^t\frac{ds}{W_{pa,i}^u(s)-W_{pa,j}^u(s)}$$
\no and $\nu_{pa}^u$-a.s. for every $t\in[0,1]$
$$W_{pa,1}^u(t)\leq...\leq W_{pa,n}^u(t)$$
\no Since $\nu_{pa}^u\sim\nu_{pa}$, $W^u$ is $\nu_{pa}$-a.s. continuous and $\nu_{pa}$-a.s. for every $1\leq i\leq n$
$$W_{pa,i}^u(t)=z_i(0)+\sigma\beta_{pa,i}^u(t)+b\int_0^tW_{pa,i}^u(s)ds+ct+\gamma\sum_{j\in\{1,...,n\}\backslash\{i\}}\int_0^t\frac{ds}{W_{pa,i}^u(s)-W_{pa,j}^u(s)}$$
\no and $\nu_{pa}$-a.s. for every $t\in[0,1]$
$$W_{pa,1}^u(t)\leq...\leq W_{pa,n}^u(t)$$

\begin{theorem} $\left(\W,\nu_{pa},\beta_{\pa},(W_{\pa}^u)_{u\in\D}\right)$ verify the
    conditions of section \ref{fr}.\newline \no  $\left(\W,\nu_{pa},\beta_{pa},\left(W_{\pa}^u\right)_{u\in G_0(\nu_{pa},\beta_{pa})}\right)$
    verify the conditions of definition \ref{DD}.
\end{theorem}

\nproof (i), (ii), (iii) and (vi) are clear. Set $u\in G_0(\nu_{pa},\beta_{pa})$, a straight calculation gives $\nu_{pa}$-a.s.
$$\beta_{pa}\circ W_{pa}^u=\beta_{pa}^u$$
hence (iv).\newline
\no Now we prove condition (v). Set $u,v\in G_0(\nu_{pa},\beta_{pa})$ such that $v+u\circ W_{pa}^v\in G_0(\nu_{pa},\beta_{pa})$, we have $\nu_{pa}$-a.s.
\beaa &&W_{pa,i}^u(t)\circ W_{pa}^v\\
&=&z_i(0)\\
&&+\left(\sigma\left(\beta_{pa,i}(t)+u(t)\right)+b\int_0^tW_{pa,i}^u(s)ds+ct+\gamma\sum_{j\in\{1,...,n\}\backslash\{i\}}\int_0^t\frac{ds}{W_{pa,i}^u(s)-W_{pa,j}^u(s)}\right)\circ W_{pa}^v\\
&=&z_i(0)+\sigma\left(\beta_{pa,i}^v(t)+u(t)\circ W_{pa}^v\right)+b\int_0^tW_{pa,i}^u(s)\circ W_{pa}^vds+ct\\
&&\indent\indent\indent+\gamma\sum_{j\in\{1,...,n\}\backslash\{i\}}\int_0^t\frac{ds}{W_{pa,i}^u(s)\circ W_{pa}^v-W_{pa,j}^u(s)\circ W_{pa}^v}\\
&=&z_i(0)+\sigma\beta_{pa,i}^{v+u\circ W_{pa}^v}(t)+b\int_0^tW_{pa,i}^u(s)\circ W_{pa}^vds+ct\\
&&+\indent\indent\indent\gamma\sum_{j\in\{1,...,n\}\backslash\{i\}}\int_0^t\frac{ds}{W_{pa,i}^u(s)\circ W_{pa}^v-W_{pa,j}^u(s)\circ W_{pa}^v}\eeaa
\no so the uniqueness of theorem \ref{par} gives $\nu_{pa}$-a.s.
$$W_{pa}^u\circ W_{pa}^v=W_{pa}^{v+u\circ W_{pa}^v}$$\nqed

\begin{corollary} It is clear that for every $u\in\D$, we clearly have $\nu_{pa}$-a.s.
$$W_{pa}^u(w)=W_{pa}^{u(w)}(w)$$
\no so theorem \ref{pl} applies.
\end{corollary}

\begin{corollary}
Theorem \ref{sde} applies. Set $u\in G_2(\nu_{pa},\beta_{pa})$, for $i\in\{1,...,n\}$, define
$$\dot{\overbrace{w_{pa,i}^u}}(t)=\dot{u}_i(t)+b\left(W_{pa,i}^u(t)-W(t)\right)+\gamma\sum_{j\in\{1,...,n\}\backslash\{i\}}\left(\frac{1}{W_{pa,i}^u(t)-W_{pa,j}^u(t)}-\frac{1}{W_i(t)-W_j(t)}\right)$$
\no We have $\nu_{pa}$-a.s.
$$W_{pa}^u=I_\W+w_{pa}^u$$
\no so
$$H(W_{pa}^u\nu_{pa}|\nu_{pa})=\frac{1}{2}\EE_{\nu_{pa}}\left[|u|_H^2\right]$$
\no if and only if there exists $v\in G_0(\nu_{pa},\beta_{pa})$ such that $W_{pa}^v$ is a strong solution to the stochastic differential system:
\beaa dW_{pa,1}^v(t)&=&\dot{\overbrace{w_{pa,1}^u}}(t)dt\circ W_{pa}^v+dW_1(t)\\
\vdots&&\\
dW_{pa,n}^v(t)&=&\dot{\overbrace{w_{pa,n}^u}}(t)dt\circ W_{pa}^v+dW_n(t)\eeaa
\end{corollary}

\vspace{2cm}
\footnotesize{
\noindent
K\'evin HARTMANN, Institut Telecom, Telecom ParisTech, LTCI CNRS D\'ept. Infres, \\
23 avenue d'Italie, 75013, Paris, France\\
kevin.hartmann@polytechnique.org}

\end{document}